\documentclass[12pt,twoside]{article}
\usepackage{amsmath}
\usepackage{amssymb}
\usepackage{enumerate,epsfig}
\newtheorem{lemma}{Lemma}[section]
\newtheorem{theorem}{Theorem}[section]

\newtheorem{defi}{Definition}[section]

\newcommand{\R}{ {\mathbb R} }

\newcommand{\ep}{{\varepsilon}}
\makeatletter
\@addtoreset{equation}{section}
\makeatother

\newcommand{\cqfd}{{\unskip\kern 6pt\penalty 500
\raise -2pt\hbox{\vrule\vbox to 6pt{\hrule width 6pt
\vfill\hrule}\vrule}\par}}
\newcommand{\eps}{{\varepsilon}}
\newcommand{\ind}{{\mathbb I}}
\begin{document}
\title{Differential Equations with singular fields}
\author{Pierre-Emmanuel Jabin\\
\footnotesize{email: jabin@unice.fr}\\
\footnotesize{Equipe Tosca, Inria, 2004 route des Lucioles, BP 93,
  06902
Sophia Antipolis}\\
\footnotesize{Laboratoire Dieudonn\'e, Univ. de Nice, Parc Valrose, 06108
  Nice cedex 02}\\}
\date{}
\maketitle
{\bf Abstract.} This paper investigates the well posedness of ordinary
differential equations and more precisely the existence (or
uniqueness) of a flow 
through explicit compactness estimates. Instead of assuming a bounded
divergence condition on the vector field, a compressibility condition
on the flow (bounded jacobian) is considered. The main result provides
existence under the condition that the vector field belongs to $BV$ in
dimension $2$ and $SBV$ in higher dimensions. 

\section{Introduction}
This article studies the existence (and secondary uniqueness) 
of a flow for the equation
\begin{equation}
\partial_t X(t,x)=b(X(t,x)),\quad X(0,x)=x. \label{ode}
\end{equation}
The most direct way to establish the existence of such of flow is of
course through a simple approximation procedure. That means taking a
regularized sequence $b_n\rightarrow b$, which enables to solve
\begin{equation}
\partial_t X_n(t,x)=b_n(X(t,x)),\quad X_n(0,x)=x,\label{sequence}
\end{equation}
by the usual Cauchy-Lipschitz Theorem. To pass to the limit in
\eqref{sequence} and obtain \eqref{ode}, it is enough to have
compactness in some strong sense (in $L^1_{loc}$ for instance) 
for the sequence $X_n$. 
Obviously some conditions are needed.

First note that we are interested in flows, which means that we are
looking for solutions $X$ which are invertible~: At least
$JX=\mbox{det}\,d_x X\neq 0\ a.e$ (with $d_x X$ the differential of
$X$ in $x$ only). So throughout this paper, only
flows $x\rightarrow X(t,x)$ which are nearly incompressible are considered
\begin{equation}
\frac{1}{C}\leq JX(t,x)\leq C,\quad \forall t\in [0,\ T],\ x\in \R^d.
\label{incomp}
\end{equation} 
If one obtains $X$ as a limit of $X_n$ then the most simple way of
satisfying \eqref{incomp} is to have
\begin{equation}
\frac{1}{C}\leq JX_n(t,x)\leq C,\quad \forall t\in [0,\ T],\ x\in
\R^d,
\label{incompn}
\end{equation}
for some constant $C$ independent of $n$. Note that both conditions
are only required on a finite and given time interval $[0,\ T]$ since
one may easily extend $X$ over $\R_+$ by the semi-group relation
$X(t+T,x)=X(t,\ X(T,x))$. Usually \eqref{incomp} and \eqref{incompn}
are obtained by assuming a bounded divergence condition on $b$ or
$b_n$ but this is not the case here.

It is certainly difficult to guess what is the optimal condition
on $b$. It is currently thought that $b\in BV(\R^d)$ is enough or (see
\cite{Br2}) 

\noindent {\bf Bressan's compactness 
conjecture~:} {\em Let $X_n$ be regular ($C^1$)
  solutions to \eqref{sequence}, satisfying \eqref{incompn} and with
  $\sup_n \int_{\R^d} |db_n(x)|\,dx<\infty$. Then the
sequence $X_n$ is locally compact in $L^1([0,\ T]\times\R^d)$.}

\medskip

From this, one would directly obtain the existence of a flow to
\eqref{ode} provided that $b\in BV(\R^d)$ and \eqref{incomp} holds.
Instead of the full Bressan's conjecture, this article essentially recovers,
through a different method, the result of \cite{ADM} namely under the
condition 
$b\in SBV$
\begin{theorem}
Assume that $b\in SBV_{loc}(\R^d)\cap L^\infty(\R^d)$ 
with a locally finite jump set (for the $d-1$ dimensional Hausdorff measure).
Let $X_n$ be regular solutions to \eqref{sequence},
satisfying \eqref{incompn} and such that $b_n\rightarrow b$ belongs uniformly to
$L^\infty(\R^d)\cap W^{1,1}(\R^d)$. Then $X_n$ is locally compact in
$L^1([0,\ T]\times\R^d)$.\label{mainresult} 
\end{theorem}
Only in dimension $2$ is it possible to be more precise
\begin{theorem}
Assume that $d\leq2$, $b\in BV_{loc}(\R^d)$. Let  $X_n$ be regular
solutions to \eqref{sequence}, 
satisfying \eqref{incompn} and such that $b_n\rightarrow b$ belongs uniformly to
$L^\infty(\R^d)\cap W^{1,1}(\R^d)$ with $\inf_K b_n\cdot B>0$
for any compact $K$ and some $B\in W^{1,\infty}(\R^d,\R^d)$. Then
$X_n$ is locally compact in 
$L^1([0,\ T]\times\R^d)$. \label{2dtheorem}
\end{theorem}
 The proof of the first result is found in section \ref{ssbv} and the
 proof of the second in section \ref{sec2d}. 
After notations and examples in section \ref{elementary} and
 technical lemmas in section \ref{technical},
particular cases are studied. In section \ref{secW11}, a very simple
proof is given if $b\in W^{1,1}$. Section \ref{sec1d} studies
\eqref{ode} in dimension 1 for which compactness holds under very
general conditions (essentially nothing for $b$ and a much weaker
version of \eqref{incomp}).  The final section offers some comments on
the unresolved issues in the full $BV$ case.

The question of uniqueness is deeply connected to the existence and in
fact the proof of Theorem \ref{mainresult} may be slighty altered in
order to provide it (it is more complicated for Th. \ref{2dtheorem}). 
Proofs are always given for the
compactness of the sequence but it is indicated and briefly explained 
after the stated
results whether they can also give uniqueness~; This is usually the
case except for sections \ref{sec1d} and \ref{sec2d}.

\bigskip 

The well posedness of \eqref{ode} is classically obtained by
the Cauchy-Lipschitz theorem. This is based on the simple estimate
\begin{equation}
|X(t,x+\delta)-X(t,x)|\leq |\delta|\,e^{t\,\|db\|_{L^\infty}}.\label{CL}
\end{equation} 
Notice that a similar bound holds if $b$ is only log-lipschitz,
leading to the important result of uniqueness for the 2d
incompressible Euler system (see for
instance \cite{Li} among many other references).

The idea in this article is to get \eqref{CL} for {\em almost all}
$x$. It is therefore greatly inspired by the recent approach developed
in \cite{CD} (see also \cite{Cr}) where the authors control the
functional
\begin{equation}
\int_{\R^d}\sup_r \int_{S^{d-1}} \log
\left(1+\frac{|X(t,x+rw)-X(t,x)|}{r}\right)\,dw\,dx.\label{functional1}
\end{equation}
This allows them to get an equivalent of Th. \ref{mainresult} provided
$b\in W^{1,p}$ with $p>1$ (and in fact $db$ in $L\log L$ would
work). Here the slighty different functional 
\begin{equation}
\sup_r \int_{\R^d}\int_{S^{d-1}} \log
\left(1+\frac{|X(t,x+rw)-X(t,x)|}{r}\right)\,dw\,dx,\label{functional2}
\end{equation}
is essentially 
considered. It gives a better condition $b\in SBV$ but has the one
drawback of not implying strong differentiability for the flow at the
limit.

Other successful approaches of course exist for \eqref{ode}. 
A most important step was achieved in \cite{DL} where the well
posedness of the flow was obtained by proving uniqueness for the
associated transport equation
\[
\partial_t u+b\cdot\nabla u=0,
\] 
under the 
conditions that $b$ be of bounded divergence and in $W^{1,1}$. The
crucial concept there is the one of renormalized solutions, namely
weak solutions $u$ s.t. $\phi(u)$ is also a solution to the same
transport equation.
This
was extended in \cite{Li2}, \cite{LL} and \cite{HLL}~; also see
\cite{BC} where the connection between well posedness for the
transport equation and \eqref{ode} is thoroughly analysed, both for
bounded divergence fields and an equivalent to \eqref{incomp}. 

Using a slight different renormalization for the equation $\partial_t
u+\nabla\cdot(b u)=0$, the well posedness  was famously obtained in \cite{Am}  
under the same bounded divergence condition and $b\in BV$. This last
assumption $b\in BV$ was also considered in \cite{CL} and
\cite{CL2}. Still using renormalized solutions, the restriction of
bounded divergence was weakened in \cite{ADM} to an assumption
equivalent to \eqref{incomp}; Unfortunately this required $b\in SBV$. 

In comparison to \cite{ADM}, Theorem \ref{mainresult} is slighty
weaker (due to the assumption of bounded jump set for ${\cal
  H}^{d-1}$), except in $2d$ where \eqref{2dtheorem} is stronger
($b\in BV$ instead of $SBV$). The main advantage of the approach
presented here is that we work directly on the differential equation,
giving for instance a very simple and direct theory for $b\in
W^{1,1}$. It also provides quantitative estimates for
$|X(t,x)-X(t,x+\delta)|$, which is connected to the regularity
of the trajectories. 

Usually using transport equations for \eqref{ode} does not directly gives
estimates like \eqref{CL}, \eqref{functional1} or
\eqref{functional2}. 
Some information of this kind can still be derived, for instance by
studying the differentiability or approximate differentiability of the
flow as in \cite{AC}, \cite{ALM}.

When an additional structure is known or assumed for $b$ (not the case
here), the conditions for
existence or uniqueness can often be loosen. The typical and most
frequent example is the hamiltonian case. For Vlasov equations for
instance, uniqueness under a $BV$ condition was obtained earlier and
more easily in \cite{Bo}~; It was even derived under slighty less than
$BV$ regularity in \cite{Ha2}. In low dimension this structure is
especially useful as illustrated in \cite{BD} where uniqueness is
obtained for Vlasov equations in dimension 2 of phase space for
continuous force terms; 
In \cite{Ha1} for a general hamiltonian
(still in 2d) with $L^p$ coefficients~; And in \cite{CCR} for a
continuous $b$ in 2d but only a bounded divergence condition instead
of a complete hamiltonian structure. In many of
those situations an estimate like \eqref{CL} or its variants
\eqref{functional1} or \eqref{functional2} is simply
false however and the flow therefore less regular than what can be proved
with more regularity.

The issue of which conditions would be optimal is still open
for the most part; Of course it should also depend on which exact property
(uniqueness or something more precise like \eqref{CL}) is looked
after. Interesting counterexamples are nevertheless known to test this
optimality, from the
early \cite{Ai}, to \cite{Br1} and \cite{DeP} which
indicates that in the general framework $BV$ indeed plays a critical role.  

Finally, there are many other ways to look for solutions to \eqref{ode}, which
are not so relevant here because they do not produce flows. A well
known example is found in 
\cite{Fi} where it is noticed 
that a one sided Lipschitz condition on $b$ is enough to get
either existence or uniqueness (depending on the side) but not
both. This is usually not enough to define a flow but can be useful to
deal with characteristics for hyperbolic problems (see \cite{Da}).

Similarly if one is interested mainly in well posedness for hyperbolic
problems, other approaches than renormalization (entropy
solutions for instance) exist. We refer to
\cite{CP} and \cite{DeL}. Those do not always yield flows, especially
where nothing is assumed on the divergence of $b$, but the
characteristics for relevant
physical solutions  are not always flows~: See for instance \cite{BJ}
in connection with sticky particles, 
\cite{PR} for the use of Filippov characteristics and \cite{PP} for
the use of an entropy condition.
%
\section{Elementary considerations\label{elementary}}
\subsection{Reduction of the problem and notations}
First if $b\in L^\infty$, $|X(t,x)-x|\leq \|b\|_{L^\infty}\,t$ and as
we look at \eqref{ode} for a finite time, we may of course reduce
ourselves to the case of a bounded domain $\Omega$ and assume that
$db$ is compactly supported in $\Omega$. 

Next note that the time dependent problem
\[
\partial_t X=b(t,X(t,x)),\quad X(0,x)=x,
\]
can be reduced to \eqref{ode} simply by adding time as a
variable. This has for first consequence to increase the dimension by
$1$, which has no importance for Th. \ref{mainresult} but could matter
in low dimension trying to use the results of sections \ref{sec1d} or
\ref{sec2d}. Second it would require that $b$ be $SBV$ in $x$ {\em
  and} $t$. 

With the right result,  it would be easy to get rid of this additional time
regularity. More precisely assuming that one has a theorem like in section
\ref{sec2d} giving for $b\in BV$ 
\begin{equation}\begin{split}
\sup_\delta \int_{\Omega} \log&
\left(1+\frac{|X(t,x+\delta)-X(t,x)|}{r}\right)\,dw\,dx\\
&\leq
C(\|b\|_\infty)\,\left(|\Omega|+\int_{\R^d\times[0,\ T]} (|\partial_t b|+|d_x
b|)\,dt\,dx\right). 
\end{split}\label{perfect}\end{equation} 
Take $\eps>0$ and change variables $X_\eps=X(\eps\,t,x)$,
$b_\eps(t,x)=\eps\,b_\eps(\eps\,t,x)$
so that
\[
\partial_t X_\eps=b_\eps(t,X_\eps),
\]
and \eqref{incomp} of course holds for $X_\eps$ on the time interval
$[0,\ T/\eps]$.  
Now applying \eqref{perfect} for $X_\eps$ and 
letting $\eps$ go to $0$, one recovers at the limit an estimate
without time derivative
\[\begin{split}
\sup_\delta \int_{\Omega} \log&
\left(1+\frac{|X(t,x+\delta)-X(t,x)|}{r}\right)\,dw\,dx\\
&\leq
C(0)\,\left(|\Omega|+\int_{\R^d\times[0,\ T]} |d_x
b|\,dt\,dx\right), 
\end{split}\]
so that $b\in L^1([0,\ T],\ BV(\R^d))$ is enough. In fact for the two
theorems concerned by this remark in this paper 
(the $W^{1,1}$ case in Th. \ref{W11}
and the $2d$ case in Th. \ref{2d}), it is easy to directly modify the
proof and obtain $b\in L^1([0,\ T],\ W^{1,1}(\R^d))$ for Th. \ref{W11}
and $b\in L^1([0,\ T],\ BV(\R))$ for Th. \ref{2d}. Unfortunately this
is of no use for theorem \ref{sbv}.

Note moreover that taking $\eps=2/\|b\|_\infty$ 
one may always assume $\|b\|_{L^\infty}\leq 2$
provided $b\in L^\infty$. 

Finally, adding one variable even in the time independent case, it is
possible to have $b_1\geq 1$. This will simplify some proofs and is
crucial for topological reasons in section \ref{sec2d}.

In summary we may work with  $b$ satisfying
\begin{equation}
M=\int_{\R^d} |db(x)|\,dx<\infty,\quad \mbox{supp}\, |db|\subset \Omega,
\quad \|b\|_{L^\infty} \leq 2,\quad b_1(x)\geq 1\  \forall x\in\Omega.
\label{assumpb}
\end{equation}

The support is denoted supp, ${\cal H}^\gamma$ denotes the $\gamma$
dimensional Hausdorff measure. $C$ will denote any universal constant
(possibly depending only on the dimension $d$ and $\Omega$) and its
value may thus change from line to line.
%
\subsection{A simple $1d$ example}
As a warm up, study the usual counterexample to Cauchy-Lipschitz in
$1d$, namely take
\[
b(x)=\sqrt{|x|}.
\]
There are several solutions to \eqref{ode} with starting point
$x<0$. All solutions are the same for some time
\begin{equation}
X(t,x)=-(t/2-\sqrt{|x|})^2,\ t\leq 2\sqrt{|x|}.\label{same}
\end{equation} 
After $t=2\sqrt{|x|}$ there are infinitely many possibilities, first
\begin{equation}
X(t,x)=(t/2-\sqrt{|x|})^2,\ t\geq 2\sqrt{|x|},\label{right}
\end{equation} 
and then for any $t_0\in[2\sqrt{|x|},\ \infty]$
\begin{equation}
X(t,x)=0\quad\mbox{for}\ 2\sqrt{|x|}\leq t\leq t_0,\qquad
X(t,x)=(t/2-t_0)^2,\ t\geq t_0.\label{good}
\end{equation} 
However among all those solutions there is only one which defines a
flow (makes $X$ invertible) and it is \eqref{right}. For the point of
view followed in this paper, this is the {\em right} one but there
could be situations where it is not the relevant solution (for
physical reasons, entropy principles...) and another should be chosen
(typically the one corresponding to $t_0=\infty$).

Note that obviously div$\,b$ is not bounded but the solution
\eqref{right} is the only one to satisfy a weaken version of
\eqref{incomp}, namely \eqref{weakcomp} (see section \ref{sec1d} where
the 1d case is studied with a result containing this
example). This shows that 
 there is a selection principle hidden in \eqref{incomp}. 
%
\subsection{Comments on the compressibility condition
  \eqref{incomp} \label{incompress}} 
Instead of \eqref{incomp}, many works rather use the condition that
the divergence of $b$ is bounded. Of course if $\mbox{div}\,b\in
L^\infty$ then \eqref{incomp} holds for any solution to \eqref{ode} so
the question is only how more general \eqref{incomp} is. Two examples
are shown to try to investigate this.

First recall that \eqref{incomp} may be reformulated in terms of the
transport equation. Namely as noticed and widely used in \cite{ADM}, 
it is equivalent to the existence of a
solution $u$ to
\[
\partial_t u+\nabla\cdot(b\,u)=0,\ \inf u(t=0)>0,\ \sup u(t=0)<\infty,
\] 
and s.t. for any $t\in [0,\ T]$
\begin{equation}
\sup_x u(t,x)\leq C\inf_x u(t=0),\quad \inf_x u(t,x)\geq C^{-1}\sup_x u(t=0).
\label{maximum}
\end{equation}
In general it it very difficult to obtain such bounds without an
assumption on the divergence. However say that $b$ is itself computed
thanks to an equation~: $b(t,x)=\nabla A(u(t,x))$ where $u$ is an
already obtained solution to the hyperbolic problem
\[
\partial_t u+\nabla\cdot (A(u))=0.
\]
Then the maximum principle for this hyperbolic equation directly gives
\eqref{maximum}. In this case \eqref{incomp} is more natural than a
bounded divergence hypothesis. It is nevertheless still usually possible to
use this latter assumption by adding the time as a dimension and
considering the stationary problem. See \cite{BC} for a full and
general analysis of the connection between transport equations and ODE's.

\bigskip

As a second example, consider what \eqref{incomp} implies for $b$
where it is discontinuous. Indeed if $\mbox{div}\,b\in L^p$ and $b$ is
discontinuous across a regular hypersurface $H$ then it is well known that the
normal component $b\cdot \nu$ ($\nu$ being the normal to $H$) {\em
  cannot} jump across $H$. This is not true anymore with only
\eqref{incomp}. Take the simple case in dimension $1$
\[
b(x)=1\ \mbox{if}\ x<0,\quad b(x)=1/2\ \mbox{if}\ x>0.
\] 
Then solving \eqref{ode} gives 
\[\begin{split}
&X(t,x)=x+t\ \mbox{if}\ t\leq -x,\quad X(t,x)=(t+x)/2\ \mbox{if}\ t\geq
-x,\\
& X(t,x)=x+t/2\ \mbox{if}\ x>0,
\end{split}\] 
which obviously satisfies \eqref{incomp}. So it can be seen that the
condition \eqref{incomp} imposes less constraint on the jumps of $b$.

Note that conversely it can be shown
that if $b\in BV$, has a discontinuity along $H$ and denoting $b^-$
and $b^+$ the two traces (see \cite{AFP}) 
then \eqref{incomp} implies that $b^+\cdot
\nu$ and $b^-\cdot \nu$ have the same sign and 
\[
b^+(x)\cdot \nu(x)\geq C^{-1}\,b^-(x)\cdot \nu(x),\quad \mbox{for}\
{\cal H}^{d-1}\ all\ x\in H.
\]
Indeed take any $x_0\in H$ s.t. $b^-$ and $b^+$ are approximately
continuous at $x_0$ (see \cite{AFP}). As $H$ is regular, change
variable so that around $x_0$, $H$ has equation $x_1=0$ (this may
modify the constant $C$ in \eqref{incomp} which hence depends on
$H$). Now consider the domain $\omega_{r,\eta}$ defined by
\[
\omega_{r,\eta}=\{x,\ -\eta<x_1<\eta,\ |x-x_0|\leq r\},
\]
and its border
\[
\omega_{r,\eta}^{\pm}=\{x,\ x_1=\pm \eta,\ |x-x_0|\leq r\},\quad
\omega_{r,\eta}^0=\{x,\ -\eta<x_1<\eta,\ |x-x_0|= r\}.
\]
Choose $r$ small enough and a sequence $\eta_n$ by approximate
continuity s.t. 
\[
\frac{1}{|\omega_{r,\eta_n}^\pm|}\int_{\omega_{r,\eta_n}^\pm} 
|b(x)-b^\pm(x_0)|\,d{\cal H}^{d-1}(x)\leq \frac{1}{C} |b^\pm_1(x_0)|.
\]
Finally compute 
\[
V(t)=\int \ind_{X(t,x)\in\omega_{r,\eta_n}}\,dx,
\]
First by simply changing variables and using \eqref{incomp} to get
\[
V(t)\leq C|\omega_{r,\eta_n}|\leq C\,r^{d-1}\,\eta_n.
\]
And then by
\[\begin{split}
&\frac{dV(t)}{dt}=\int \delta(X(t,x)\in\omega_{r,\eta_n}^0)\,
b(X(t,x))\cdot (X'-x_0)\,r^{-1}\\
&\quad+ \int \delta(X(t,x)\in\omega_{r,\eta_n}^+)\,
b_1(X(t,x))\\
&\quad- \int \delta(X(t,x)\in\omega_{r,\eta_n}^-)\,
b_1(X(t,x)),\\
\end{split}\]
where $X'=(0,X_2,\ldots,X_d)$.
The first term is bounded by
$
\|b\|_\infty \eta_n.
$
Assume for instance that $b_1^+(x_0)>0$ then by change of variables
the second term is smaller than $C\,{\cal H}^{d-1}(\omega_{r,\eta_n}^+)
\,b^+_1(x_0)=C\,r^{d-1} \,b^+_1(x_0)$. The third term is larger than
$r^{d-1} b^-_1(x_0)/C$ (if $b^-_1(x)>0$). Integrating over $[0,\ T]$
and taking $\eta_n$ small enough leads to
\[
b^+_1(x_0)\geq C^{-1}\,b^-_1(x).
\]

\section{Preliminary results \label{technical}}
Two simple lemmas are given here, which will be used frequently in
other proofs. 
\begin{lemma} Assume $b\in BV$. 
There exists a constant $C$ s.t. for any $x$, $y$
\begin{equation}
|b(x)-b(y)|\leq C\,\int_{B(x,y)}
|db(z)|\,\left(\frac{1}{|x-z|^{d-1}}+\frac{1}{|y-z|^{d-1}}
\right)\,dz,\label{diff}
\end{equation}
where $B(x,y)$ denotes the ball of center $(x+y)/2$ and diameter
$|x-y|$. \label{lemma}
\end{lemma}
\noindent{\bf Proof.}
This is just an  explicit computation: Change
coordinates so that $x=(-\alpha,0,\ldots)$ and $y=-x$. Then take a
path $t\in [0,\ 1/2]\longrightarrow (-\alpha+2\alpha t, \alpha t r)$
for any $r$ in the unit ball of $\R^{n-1}$ and take the symmetric path
for $t>1/2$. Then all those paths $\gamma_r$ connect $x$ and $y$
so that
\[
|b(x)-b(y)|\leq \int_{\gamma_r} |db(z)|\,dl(z).
\]
Averaging
over $r$ in the ball $B(x,y)$ and changing coordinates get
\eqref{diff}.\cqfd

A slight variant of this (usefull for the $SBV$ case in particular) is
\begin{lemma} Assume $b\in BV$ and $H$ is an hypersurface, lipschitz
  regular of $\R^d$. 
There exists a constant $C$ and a constant $K$ (depending on $H$) 
s.t. for any $x$, $y$ locally on the same side of $H$
\begin{equation}
|b(x)-b(y)|\leq C\,\int_{B_K(x,y)\setminus H}
|db(z)|\,\left(\frac{1}{|x-z|^{d-1}}+\frac{1}{|y-z|^{d-1}}
\right)\,dz,\label{diff2}
\end{equation}
where $B_K(x,y)$ denotes the ball of center $(x+y)/2$ and diameter
$K\,|x-y|$. \label{lemma2}
\end{lemma}
\noindent{\bf Proof.} This is just the same idea as before~: Consider
all paths $\gamma$ connecting $x$ and $y$, of length at most
$K\,|x-y|$ and not crossing $H$. Average over those to get the result.

Note that $K$ must be larger than the lipschitz regularity of $H$: If
$H$ is locally given by the equation $f(x)=0$ then $K\geq
C\,\|df\|_{\infty}$. And if $K$ is chosen like that then the
definition of locally on the same side can simply be: There exists a
path $\gamma$ connecting $x$ and $y$ without crossing $H$ and of
length less than $3K/2\,|x-y|$.\cqfd

\bigskip

Finally let us note that Lemma \ref{lemma} is a more precise version
of the well known bound (used in \cite{CD} in particular)
\begin{equation}
|b(x)-b(y)|\leq C\,|x-y|\,(M|db|(x)+M|db|(y)),\label{maximal}
\end{equation}
where $M|db|$ is the maximal function of $|db|$. Indeed decomposing
$B(x,y)$ into $\bigcup (R(x, 2^{-n})\cap B(x,y))$ for all $n\geq
n_0=-\log_2 |x-y|$, with $R(x,r)=\{z,\ r/2\leq |z-x|< r\}$, one gets
\[\begin{split}
\int_{B(x,y)} \frac{|db(z)|}{|x-z|^{d-1}}\,dz&\leq \sum_{n\geq n_0}
2^{(d-1)(n+1)} \int_{R(x,2^{-n})} |db(z)|\,dz\\
&\leq 2^{d-1} \sum_{n\geq
  n_0} 2^{-n} M|db(x)|\,dz\leq 2^{d}\,|x-y|\,M|db|(x),
\end{split}\]
recalling that
\[
M|db(x)|=\sup_r r^{d}\int_{B(x,r)} |db(z)|\,dz.
\]
%
\section{The $W^{1,1}$ case\label{secW11}}
\subsection{The result}
Following  \cite{CD}, we define for any $\delta$
\[
Q_\delta(t)=\int_{\Omega} 
\log\left(1+\frac{|X(t,x)-X(t,x+\delta)|}{|\delta|}\right)\,dx.
\]
As $db$ belongs to $L^1$, there exists $\phi\in C^\infty(\R_+)$ with
\begin{equation}
  \phi(\xi)/\xi\ \mbox{increasing},\quad
  \frac{\phi(\xi)}{\xi}\longrightarrow +\infty\ as\ \xi\rightarrow +\infty,
\label{propphi}
\end{equation}
and such that
\begin{equation}
\int_{\Omega} \phi(|db(x)|)\,dx<\infty.\label{bphi}
\end{equation}
The main result here is the explicit estimate
\begin{theorem}
Assume that $b$ satisfies \eqref{incomp} and \eqref{bphi} then
there exists a constant $C$ depending only on $\Omega$  and a
continuous function
$\psi$ depending only on $\phi$ and with {\bf $\psi(\xi)/|\log
  \xi|\rightarrow 0$} as $\xi\rightarrow 0$ such that
\[
Q_\delta(t)\leq |\Omega|\,\log 2+C\,t\,\psi(|\delta|)\,
\int_{\Omega} (1+\phi(|db(x)|))\,dx.
\]\label{W11}
\end{theorem}
\noindent {\bf Remarks.} 

1. This of course immediately implies that any sequence of 
solutions $X_n$ to \eqref{sequence} is compact thus proving Bressan's
conjecture in the restricted $W^{1,1}$ case.\\
2. Uniqueness~: The proof is identical if one considers two different
solutions $X$ 
and $Y$ to \eqref{ode}. Therefore the solution is also unique.

\medskip

\noindent{\bf Proof of Th. \ref{W11}.} \\
Start by differentiating $Q_\delta$ in time
\[
Q_\delta'(t)\leq \int_{\Omega} \frac{|\partial_t X(t,x)-\partial_t
  X(t,x+\delta)|}{|\delta|+|X-X_\delta|}\,dx,
\]
where $X_\delta$ stands for $X(t,x+\delta)$. As $X$ solves
\eqref{ode},
\[
Q_\delta'(t)\leq \int_{\Omega} \frac{|b( X(t,x))-
  b(X(t,x+\delta))|}{|\delta|+|X-X_\delta|}\,dx,
\]
and using lemma \ref{lemma}, one gets
\[
Q_\delta'(t)\leq \int_{x,z} \ind_{z\in B(X,X_\delta)}\,
\frac{|db(z)|}{|\delta|+|X-X_\delta|}\,
\left(\frac{1}{|z-X|^{d-1}}+\frac{1}{|z-X_\delta|^{d-1}}\right)\,dx\,dz. 
\]
Now for any $M$ decompose the integral in $z$ into the domain $E_M$
where $|db(z)|\leq M$ and the set $F_M$ where $|db(z)|\geq M$
\[\begin{split}
Q_\delta'(t)\leq I+II,
\end{split}\] 
with
\[\begin{split}
I&=\int_{\Omega}\!\!\int_{E_M} \ind_{z\in B(X,X_\delta)}\,
\frac{|db(z)|}{|\delta|+|X-X_\delta|}\,
\left(\frac{1}{|z-X|^{d-1}}+\frac{1}{|z-X_\delta|^{d-1}}\right)\,dx\,dz\\
&\leq \int_{x,z}\ind_{z\in B(X,X_\delta)}\,
\frac{M}{|\delta|+|X-X_\delta|}\,
\left(\frac{1}{|z-X|^{d-1}}+\frac{1}{|z-X_\delta|^{d-1}}\right)\,dx\,dz
\\
&\leq M |\Omega|.
\end{split}\]
On the other hand
\[\begin{split}
II&=\int_{\Omega}\!\!\int_{F_M} \ind_{z\in B(X,X_\delta)}\,
\frac{|db(z)|}{|\delta|+|X-X_\delta|}\,
\left(\frac{1}{|z-X|^{d-1}}+\frac{1}{|z-X_\delta|^{d-1}}\right)\,dx\,dz\\
&\leq \frac{M}{\phi(M)}\,\int_{x,z}
\phi(|db(z)|)\,\Biggl(\frac{1}{(|\delta|+|z-X|)\,(|z-X|^{d-1})}\\
&\qquad\qquad\qquad +\frac{1}{(|\delta|+|z-X_\delta|)
\,(|z-X_\delta|^{d-1})}\Biggr)\,dx\,dz,
\end{split}\]
since first as $\phi/\xi$ is increasing, if $|db|>M$ then $|db|\leq
\phi(|db|)\, M/\phi(M)$ and second as $z\in B(X,X_\delta)$ then
$|z-X|\leq |X-X_\delta|$ and $|z-X_\delta|\leq |X-X_\delta|$.

Note that $X$ and $X_\delta$ play the same role and in particular the
transform $x\longrightarrow X_\delta=X(t,x+\delta)$ also has a bounded
jacobian. Therefore we change variable in $x$, for the first term in
the parenthesis from $x$ to $X$ and for the second from $x$ to
$X_\delta$ to find
\[\begin{split}
II&\leq 2\,\frac{M}{\phi(M)}\,\int_{x,z}
\phi(|db(z)|)\,\frac{1}{(|\delta|+|z-x|)\,(|z-x|^{d-1})}
\,dx\,dz\\
&\leq 2\,\frac{M}{\phi(M)}\,\log(1/|\delta|)\,\int \phi(|db(z)|)\,dz, 
\end{split}\]
by integrating first in $x$. Combining both estimates, one obtains
\[
Q_\delta'(t)\leq (M+2\,\frac{M}{\phi(M)}\,\log(|\delta|^{-1}))
\,\int (1+\phi(|db(z)|))\,dz.
\]
Defining $\psi(|\delta|)=\inf_M
M+2\,\frac{M}{\phi(M)}\,\log(|\delta|^{-1})$, this concludes the
proof.\cqfd

\medskip

Note that this proof uses a sort of interpolation of $db$ between
$L^\infty$ and $L^1$. If instead one uses the result of \cite{CD},
then it is enough to interpolate between $L\log L$ and $L^1$ and $\psi$
is then defined by
\[
\psi(|\delta|)=\inf_M \frac{M\log M}{\phi(M)}
+2\frac{M}{\phi(M)}\,\log(|\delta|^{-1}),
\]
provided that $\phi(\xi)\leq \xi\,log\xi$. This last estimate for
$\psi$ is of course much better than the previous one, even though for
well posedness (compactness or uniqueness), it does not matter.
\subsection{An example}
The following remark was first made by S. Bianchini. The
previous proof shows that if $b\in W^{1,1}$ then for any flow $X$
satisfying \eqref{incomp}, one may bound
\begin{equation}
\int_\Omega \frac{|b(X(x))-b(X(x+\delta))|}{|\delta|+|X-X(x+\delta)|}\,dx,
\label{linear}
\end{equation}
by $o(\log |\delta|)$. It is then essentially a linear estimate in the
sense that in \eqref{linear} one never uses that $b$ and $X$ are
connected. The situation below shows that this simple way of controlling
\eqref{linear} cannot be extended further to $b\in BV$.

The example is shown in $1d$ but it can of course be extended to any
dimension. Simply take for $b$ the Heaviside step function
$\ind_{x>0}$. As for $X$, fix $n$ and choose
\[\begin{split}
&X(t,x)=x,\quad if\ x\in [k/n,\ (2k+1)/2n]\\
&\qquad\qquad\mbox{or}\ x\in [-(2k+1)/2n,\
-k/n],\ 1\leq k\leq n-1\\
& X(t,x)=-x,\quad if\ x\in [(2k+1)/2n,\ (k+1)/n]\\
&\qquad\qquad \mbox{or}\quad x\in[-(k+1)/n,\
-(2k+1)/2n],\ 1\leq k\leq n-1,\\
& X(t,x)=x\ \mbox{otherwise}.
  \end{split}
\]
It is obvious that $X$ satisfies \eqref{incomp} (it only swaps the
intervals) and choosing $\delta=1/2n$
\[
\int_{-1}^1
\frac{|b(X(x))-b(X(x+\delta))|}{|\delta|+|X-X(x+\delta)|}\,dx\geq
C \sum_{k=1}^{n-1}\; \frac{1}{2n} \frac{n}{k}\geq C\log n=C\log |\delta|.
\]
To bypass this obstacle, it is necessary to use the dependance of $X$ in
terms of $b$ or the fact that $X(0,x)=x$. In dimension $1$ for
example, it is not possible to pass continuously from $X=x$ at $t=0$
to an $X$ as shown here.
\section{The 1d case\label{sec1d}}
The stationary equation \eqref{ode} in only one dimension is of course
a very particular situation. In this case assumption
\eqref{incomp} is more than enough and no additional regularity is
required (just as the divergence controls the whole gradient). In fact
it is even too much and one only needs to assume that the image of a non
empty interval remains non empty
\begin{defi}
{\bf Weak compressibility~:} $\exists \phi\in C(\R_+)$ with
$\phi(\xi)>0$ for any $\xi>0$ s.t. $\forall I$ interval, $\forall
t\in[-T,\ T]$
\[
|X(t,I)|>\phi(I).
\]  
\label{weakcomp}
\end{defi}
Then it is quite straightforward to get
\begin{theorem}
There exists a strictly increasing $\psi$ with $\psi(0)=0$ s.t. 
any $X$ limit of 
solutions to \eqref{sequence} with $b_n\in W^{1,\infty}_{loc}$ (but
not necessarily the limit $b$), and satisfying \eqref{weakcomp} also satisfies
\[
|X(t,x+\delta)-X(t,x)|\leq \phi(\delta).
\]
\end{theorem}
{\bf Remarks.}\\
1. This is of course enough to ensure compactness and Bressan's
conjecture in this case.\\
2. The proof relies on the topology of $\R$. It does not give directly
any
uniqueness result and in particular, $X(t,x+\delta)$ cannot be
replaced by another solution $Y$. Of course once the regularity of the
solution is obtained it should be possible to then derive uniqueness.

\medskip

{\bf Proof.} It is enough to prove that the estimate holds for a
regular solution of \eqref{ode}. First note that as if $I\subset J$,
$X(t,I)\subset X(t,J)$ then we may always take $\phi$ increasing
in \eqref{weakcomp}. 

As the solution to \eqref{ode} is unique, the image of the
interval $[x_1,\ x_2]$ is simply the interval $[X(t,x_1),\
X(t,x_2)]$. 
This also means that the image by the flow $X(t,.)$ of
the interval $[X(t,x),\
X(t,x+\delta)]$ is the interval $[x,\ x+\delta]$ so applying
\eqref{weakcomp}
\[
\delta>\phi(|X(t,x+\delta)-X(t,x)|),
\]
which gives the result after composition with $\phi^{-1}$.\cqfd

\section{The 2-d case\label{sec2d}}
The situation in two dimension is more complicated than in one
dimension but still very constrained by the topology. 

Let us again consider the functional
\[
Q_\delta(t)=\int_{\Omega} 
\log\left(1+\frac{|X(t,x)-X(t,x+\delta)|}{|\delta|}\right)\,dy\,dx,
\]
where $\delta$ is a fixed vector, for example $\delta=(0,r)$.

In this setting it is possible to obtain the optimal
\begin{theorem}
Let $X$ be a regular solution to \eqref{ode} satisfying \eqref{incomp}
and assume that $b$ satisfies 
\eqref{assumpb}  then there exists a constant $C$
(depending only on the constant in \eqref{incomp}) such that
\begin{equation}
Q_\delta(t)\leq |\Omega|\,\log 2+C\,(t+|\delta|)\;\int_{\R^d} |db(x)|\,dx. 
\label{2d}\end{equation}\label{2dresult}
\end{theorem}
{\bf Remarks.}\\
1. The proof uses some ideas developped together with U. Stefanelli and
C. DeLellis.\\
2. The assumption $b_1\geq 1$ is {\em crucial} and deeply connected to
the $2d$ topology. Of course the constants $1$ and $2$ in
\eqref{assumpb} can be changed to $b_1\geq c_1$ and $|b|\leq c_2$ but
then the constant in \eqref{2d} is modified by $c_2/c_1$ (just a
scaling argument). Finally as
it is always possible to decompose $\Omega$ into smaller domains, this
assumption can be optimized to assuming that there are a finite number
of regular domains $\Omega_i$ with $\Omega=\bigcup \Omega_i$ and in
each $\Omega_i$ either $b\cdot B>0$ for a constant $B$ which in turn
would give the assumption in Th. \ref{2dtheorem}.\\
3. This result is optimal as the estimate in \eqref{2d} does not
depend at all on $\delta$. I have no idea how to extend it in higher
dimensions even for $b\in W^{1,1}$.\\
4. As in the $1d$ case this implies directly Bressan's compactness
conjecture but not uniqueness of the flow. This uniqueness should
be deduced in a second step using this estimate.
  
\noindent{\bf Proof of Theorem \eqref{2dresult}.}
Compute
\[\begin{split}
\frac{d}{dt} &\log\left(\frac{|X(t,x)-X(t,x+\delta)|}{|\delta|}\right)\leq 
\frac{|b(X(t,x))-b(X(t,x+\delta))|}{|\delta|+|X(t,x)-X(t,x+\delta)|}\\
&\leq \frac{|b(X(t_\delta(t,x),x))-b(X(t,x+\delta))|}{|\delta|}+
\frac{|b(X(t,x))-b(X(t_\delta(x),x))|}{|\delta|+|X(t,x)-X(t,x+\delta)|},
\end{split}\]
where $t_\delta(t,x)$ is defined as the unique time such that
\[
X_1(t_\delta(t,x),x)=X_1(t,x+\delta).
\]
Integrating over $\Omega$ and $[0,\ t]$, we find that
\[
Q_\delta(t)\leq I+II,
\]
with
\begin{equation}
I=\int_\Omega\int_0^t
\frac{|b(X(t_\delta(s,x),x))-b(X(s,x+\delta))|}{|\delta|}\,ds\,dx , 
\label{I}
\end{equation}
and
\begin{equation}
II=\int_\Omega\int_0^t
\frac{|b(X(s,x))-b(X(t_\delta(s,x),x))|}{|\delta|+|X(s,x)-X(s,x+\delta)|}\,
ds\,dx. 
\label{II}
\end{equation}

The first term may be bounded by

\[
I\leq \frac{1}{|\delta|}\,
\int_\Omega\int_0^t \int_{X_2(t_\delta,x)}^{X_2(s,x+\delta)} 
|db(X_1(s,x+\delta),\alpha)|\,d\alpha\,ds\,dx.
\]
Denote by $\Omega_\delta(x)$ the set of points included between the two
trajectories $X(s,x)$ and $X(s',x+\delta)$ for $-t\leq s,s'\leq t$ (note
that as $b_1>0$ each trajectory is indeed a 1-d manifold). 
As $b_1\geq 1$ then
$\partial_t X_1(t,x)\geq 1$ and a line of equation $x_1=\alpha$ may
cross the trajectory $\{X(s,x), s\in \R)\}$ only once. Therefore we
may describe $\Omega_\delta$ like
\[
\Omega_\delta(x)=\{(y_1,y_2)
\;|\ X_2(t(y_1),x)<y_2<X_2(t_\delta(t(y_1),x),x+\delta)\},
\]
where $t(\alpha)$ is defined as the unique $t$ such that
$X_1(t(\alpha),x)=\alpha$.

Consequently, changing variables, we get that
\[
I\leq \frac{2}{|\delta|}\int_\Omega\int_{\Omega_\delta(x)} |db(y)|\,dy\,dx,
\]
as the Jacobian of the transform $t\rightarrow X_1(t,x)$ is at most
$2$.

Changing the order of integration, we have
\[
I\leq \frac{2}{|\delta|}\int_{\tilde \Omega} |db(y)| \times|\{x,\ y\in
\Omega_\delta(x)\}|\,dy, 
\]
with $\tilde \Omega=\cup_x \Omega_\delta(x)$ and therefore $|\tilde
\Omega|\leq C |\Omega|$ (if $\Omega$ is regular).

On the other hand
\[
\Omega_\delta(x)\cap \Omega_\delta((x_1+\alpha,x_2+\beta))=\emptyset,
\]
if $|\alpha|>4t+r$ or $|\beta|>r$. This is one point
where the two dimensional aspect is crucial.

Consequently
\[
|\{x,\ y\in
\Omega_\delta(x)\}|\leq 2r\,(4t+r),
\]
and 
\[
I\leq 4\,(4t+r)\int_{\tilde \Omega} |db(y)|\,dy.
\]

\medskip

Let us now bound $II$. We first obtain
\[
II\leq \int_\Omega \int_0^t \frac{2}{|\delta|+|X(s,x)-X(s,x+\delta)|}
\left|\int_{s}^{t_\delta(s,x)}
|db(X(u,x))|\,du
\right|\,ds\,dx.
\]
Next note that
\[\begin{split}
|X(s,x)-X(s,x+\delta)|&\geq
 |X_1(s,x)-X_1(s,x+\delta)|\\&=|X_1(s,x)-X_1(t_\delta(s,x),x)|
\geq 
\frac{1}{2} |s-t_\delta(s,x)|, 
\end{split}\]
so that
\[
II\leq \int_\Omega \int_0^t \frac{2}{|\delta|+|s-t_\delta(s,x)|}
\left|\int_{s}^{t_\delta(s,x)}
|db(X(u,x))|\,du
\right|\,ds\,dx.
\]
By Fubini's Theorem, we get
\[
II\leq \int_\Omega \int_0^{t+{|\delta|}} |db(X(u,x))|
\int_0^t
\frac{\ind_{ u\in [s,\ t_\delta(s,x)]}}{r+|s-t_\delta(s,x)|}\,ds\,du\,dx.
\]
Note that the convention $[a,b]=[b,a]$ if $a>b$ is used.

First remark that, due to \eqref{assumpb}
\[
\partial_t (X_1(t_\delta(t,x),x))=\partial_t X_1(t,x+\delta)\in\
	[1,\ 2].
\]
As such
\[
b_1(\Phi(t_\delta(t,x),x))\times \partial_t t_\delta(t,x)\in\ [1,\ 2],
\]
and thanks to \eqref{assumpb} again
\begin{equation}
\partial_t t_\delta(t,x)\in\ [1/2,\ 2].\label{derivativet}
\end{equation}
Hence we define as $s_\delta(u,x)$ the unique $s$ such that
\[
t_\delta(s,x)=u.
\]
Assume that $s_\delta(u,x)\leq u$ (the other case is dealt with in the
same manner). Then $u\in [s,\ t_\delta(s,x)]$ iff $s\in [s_\delta(u,x),\
  u]$.

As long as $u\in [s,\ t_\delta(s,x)]$, we have that
\[
|s-t_\delta(s,x)|=|s-u|+|t_\delta(s,x)-u|\geq \max(|s-u|,|t_\delta(s,x)-u|).
\]
Moreover
\[
|s-u|\geq |s_\delta(u,x)-u|-|s-s_\delta(u,x)|,
\] 
and using \eqref{derivativet}
\[
|t_\delta(s,x)-u|=|t_\delta(s,x)-t_\delta(s_\delta(u,x),x)|
\geq \frac{1}{2} |s-s_\delta(u,x)|.
\]
So we bound from below, using $|s-u|$ if $|t_\delta(s,x)-u|\leq
|s_\delta(u,x)-u|/3$ 
and $|t_\delta(s,x)-u|$ otherwise
\[
|s-t_\delta(s,x)|\geq \frac{1}{3} |s_\delta(u,x)-u|.
\]
Finally this gives that
\[
\int_0^t
\frac{\ind_{ u\in [s,\ t_\delta(s,x)]}}{r+|s-t_\delta(s,x)|}\,ds\leq
3 \int_{s_\delta(u,x)}^u \frac{ds}{|s_\delta(u,x)-u|}\leq 3. 
\]
And coming back to $II$ and using \eqref{incomp}
\[
II\leq 3\int_\Omega \int_0^{t+|\delta|} |db(X(u,x))|\,du\,dx\leq 3C(t+|\delta|)
\int_{\R^d} |db(y)|\,dy,
\]
which, summing with $I$, exactly gives the theorem.\cqfd
\section{The $SBV$ case \label{ssbv}}
\subsection{Presentation}
%
We recall the definition of $SBV$ (see for example \cite{AFP}, {\em 4.1})
\begin{defi}
$b\in SBV(\R^d)$ iff $db=m+\theta {\cal H}^{d-1}|_{J}$ with $m\in
L^1$, $J$
$\sigma-finite$ with respect to ${\cal H}^{d-1}$ and
\[
\int_J |\theta|\; d{\cal H}^{d-1} <\infty. 
\]
\end{defi}
For this restricted class of $b$ but now in any dimension, Bressan's
compactness conjecture holds and more precisely
\begin{theorem} Consider a sequence of solutions
to \eqref{sequence} satisfying \eqref{incompn}, \eqref{assumpb}
uniformly in $n$ and such that $b_n\longrightarrow b\in
SBV(\R^d)$ with ${\cal H}^{d-1}(J)<\infty$. 
Then this sequence is compact and more precisely
$\forall \eta$, there exists a continuous function
$\ep(\delta)$ with $\ep(0)=0$ and
such that $\forall n$ large enough,
$\forall \delta'<\delta$, $\forall w\in S^{d-1}$,  
$\exists \omega$ (depending on $\eta$, $n$ and $\delta'$) 
with $|\omega|\leq \eta$ and 
\begin{equation}
\forall x\in\Omega\setminus\omega,\quad \forall t\leq 1,\qquad
|X_n(t,x)-X_n(t,x+\delta'w)|\leq \ep(\delta).
\end{equation} \label{sbv}
\end{theorem}
\noindent{\bf Remarks.}\\
1. Uniqueness also holds, the proof being the same. It is even easier
as there is no need to work with a fixed scale $\delta'$ and the
additional assumption ${\cal H}^{d-1}(J)<\infty$ is not required (see
the more detailed comment below).\\
2. The function $\ep(\delta)$ strongly depends on the structure of $b$
and in particular on the local regularity of its jump set $J$ (more
precisely the lipschitz norm of $g$ if $J$ has equation $g(x)=0$
locally). Therefore this result cannot be extended directly to get to
$b\in BV$.

\bigskip

Let us comment more on the assumption ${\cal 
  H}^{d-1}(J)<\infty$. A natural idea to try to bypass it would be to
truncate $J$ 
into a set with finite Hausdorff measure and a remainder $J'$. This
means that we are approximating $b$ by $b_\gamma$ with
$|b-b_\gamma|<\gamma$ and the jump set of $b_\gamma$ is a nice $J_\gamma$
with ${\cal 
  H}^{d-1}(J_\gamma)<\infty$ (and in fact it is even $o(\gamma^{-1})$).
To give
an idea of why this is not directly working, 
consider step 3 in the proof. Its aim is to 
control the number of times a trajectory $X_n(t,x)$ comes at a
distance $\delta$ of $J_\gamma$. For that the crucial estimate is
\[
|\{x,\ d(x,J_\gamma)<\delta\}|\leq K\,\delta\,{\cal H}^{d-1}(J_\gamma).
\]
However the constant $K$ in this estimate depends on $J_\gamma$. What
is true is that $|\{x,\ d(x,J_\gamma)<\delta\}|/\delta$ is bounded,
{\em asymptotically} as $\delta\rightarrow 0$,  by ${\cal
  H}^{d-1}(J_\gamma)$. But if $\delta$ is not small enough, then the
constant $K$ can be much larger than $1$~: Take $J_\gamma$ composed of
many small pieces of radius much smaller than $\delta$ for example. 

So
the only solution is to take $\delta$ small enough for
$J_\gamma$. Unfortunately we approximated $b$ so in any case we cannot
take scales smaller than $\gamma$ and of course it could very well be
that $\delta=\gamma$ is still not small enough...

On the other hand this is a problem only when considering a positive
scale. 
If the aim is only uniqueness, one does
not have to choose a scale and following the same steps, it is enough
to control the number of times that $X(t,x)$ crosses $J_\gamma$. This
number is always bounded directly in terms of ${\cal
  H}^{d-1}(J_\gamma)$. This would make the proof of uniqueness really
easier and without the assumption ${\cal H}^{d-1}(J)<\infty$, in line
with \cite{ADM}.

\bigskip

Finally, before giving the details, let us present the main ideas of
the proof. The contributions from the jump part of $db_n$ and from the
$L^1$ part will be treated separately and for the $L^1$ part of course
similarly
as the $W^{1,1}$ case (see \ref{secW11}). So focus here only on the
jump part and simply assume that $b_n$ is piecewise constant~:
$b_n(x)=b^-$ for $x_1<0$ and $b_n(x)=b^+$ for $x_1>0$ with $b^-\neq
b^+$.

Take the two trajectories $X_n(t,x)$ and $X_n(t,x+\delta w)$. And
assume that initially $x_1<-\delta$ (if $x_1>\delta$ then they never
see the jump in $b_n$). Until one of the two $X_n$ or
$X_{n,\delta}=X_n(t,x+\delta w)$ reaches the hyperplane $x_1=0$, their
velocity is the same. Assume that $X_{n,\delta}$ touches the
hyperplane first and denote $t_1$ the first time in $[0,\ T]$ when
this happens. So
\[
\forall t\leq t_1,\quad |X_n-X_{n,\delta}|=\delta.
\]
As $b_1\geq 1$ and in particular $b^+\geq 1$, $X_{n,\delta}$ will
never again pass through $\{x_1=0\}$ so that
\[
b_n(X_{n,\delta})=b^+,\ \forall t\geq t_1.
\]
However $b^-\geq 1$ also so $X_n$ will necessarily touch the
hyperplane some time after $t_1$. Denote $t_2$ this time. After $t\geq
t_2$, $b_n(X_n)=b^+$ and so
\[
|X_n(t)-X_{n,\delta}(t)|=|X_n(t_2)-X_{n,\delta}(t_2)|.
\]
As $\|b\|_\infty\leq 2$, this implies that for any $t\in [0,\ T]$
\[
|X_n(t)-X_{n,\delta}(t)|\leq \delta+2(t_2-t_1).
\]
Finally at $t_1$, $X_{n,1}(t_1)\geq -\delta$ and as $b^-\geq 1$ this
means that $t_2\leq t_1+\delta$, enabling us to conclude that
\[
|X_n(t)-X_{n,\delta}(t)|\leq \delta+2\delta,\ \forall\,t\in[0,\ T].
\]
Notice that this is only one of two cases~: Here the trajectories
always cross the jump set but it could happen that they are
tangent. So another interesting example occurs when~: $b_n(x)=b^-$ for
$x_2<0$, $b_n(x)=b^+$ for $x_2=0$ and $b^-_2=0$ (and therefore
$b^+_2=0$ by the compressibility condition \eqref{incompn}, see also
section \ref{incompress}). Now the trajectories never cross $\{x_2=0\}$ so if
$x$ and $x+\delta w$ are on the side of this hyperplane, $X_n$ and
$X_{n,\delta}$ stay on the same side and hence
\[
|X_n-X_{n,\delta}|=\delta.
\]
We have problems when they start on different sides and then there is
nothing one can do~: At time $t$, $|X_n-X_{n,\delta}|$ is of order
$t$. However this only happens if $x$ belongs to $\omega=\{-\delta\leq
x_2\leq\delta\}$, which is why in the theorem we need to exclude some
starting points. Here we would simply have
\[
\exists \omega\subset\Omega\ \mbox{with}\ |\omega|\leq C\,\delta,\quad
\forall x\in \Omega\setminus\omega,\ \forall t\in [0,\ T],\
|X_n-X_{n,\delta}|\leq \delta.  
\]

\smallskip

In the general $SBV$ case, the two situations may occur. So it is
necessary to first identify the jump set, the regions where the
trajectories typically cross and the regions where they are almost
tangent (and all that quantitatively). Next one has to exclude the
starting points $x$ which would lead to trajectories passing through
the tangent regions, and exclude among the other trajectories the ones
that are not typical ({\it i.e.} they do not cross as fast as they
should). 
And of course a
trajectory could very well cross the jump set several times so a bound
on that number of times will be required as well. 
%
\subsection{Proof of Theorem \ref{sbv}}

\noindent{\em Step 1~: Decomposition of $db$, $db_n$.}

Fix $\eta$. Through all proof $K$ will denote constants depending on
$\eta$ or $b$  and $C$ will be kept for constants depending only on $d$ or
$\Omega$.

Decompose $b$ as in the definition
\[
db=m+\theta{\cal H}^{d-1}|J.
\]
$J$ is countably rectifiable and of finite measure. 
Therefore decompose $J=H\cup J'$ with
$H$ a finite union of rectifiable sets $H_i$ such that
\begin{equation}
{\cal H}^{d-1}(J')< \eta/C.
\end{equation}
By the definition of the Hausdorff measure (and as $J'$ is countably
rectifiable), there exists a covering $J'\subset \bigcup_i
B(x_i,r_i)$ s.t. a point of $\R^d$ belongs to at most $C$ balls and with 
\[
\sum_i r_i^{d-1}<2\eta/C.
\]
As $b_n$ converges toward $b$, and taking $n\geq N$ large enough, one
may decompose accordingly $db_n$ as
\[
|db_n|=m_n+\sigma_n+r_n,\quad m_n,\ \sigma_n,\ r_n\geq 0,
\] 
with for some $\tilde \phi$ with $\tilde \phi(\xi)/\xi\rightarrow +\infty$ as
$\xi\rightarrow +\infty$
\begin{equation}
\sup_n\,\int_\Omega \tilde \phi(m_n)\,dx\leq C,\label{mn}
\end{equation}
and
\begin{equation}
\mbox{supp}\, \sigma_n\subset \{x,\ d(x,H)\leq \delta\},
\label{sigma}
\end{equation}
with finally
\begin{equation}
\mbox{supp}\, r_n\subset \bigcup_i B(x_i,2r_i).
\end{equation}
Denote $\tilde b_n$ the function equal to $b_n$ on
$\Omega\setminus\Omega_r$ with
$\Omega_r=\bigcup_i B(x_i,4r_i)$ and in $B(x_i,3r_i)$
\[
\tilde b_n=b_n\star L_{r_i},\ \mbox{with}\ L_r(x)=r^{-d} L(x/r)
\]
with $L$ a $C^\infty$ function with total mass $1$ and compactly
supported in $B(0,1)$. In $B(x_i,4r_i)\setminus B(x_i,3r_i)$, 
choose a linear interpolation between the two values on $\partial
B(x_i,4r_i)$ and $\partial B(x_i,3r_i)$. One obtains a corresponding 
decomposition of
$|d\tilde b_n|$
\[
|db_n|\leq m_n \ind_{\Omega_r^c}+\sigma_n\ind_{\Omega_r^c}
+\tilde r_n+\frac{\mu_n}{r_i},
\]
with $\mu_n\leq 2$ a bounded function and
\[
\tilde r_n\leq \sum_i |db_n|\star L_{r_i}\,\ind_{B(x_i,3r_i)}.
\]
By De La Vall\'ee Poussin, 
there exists $\Phi$, with $\psi(\xi)=\Phi(\xi)/\xi$ increasing and
converging to $+\infty$ s.t.
\[\begin{split}
\int_\Omega \Phi(m_n+\tilde r_n+\mu_n/r_i)\,dx\leq &\int_\Omega
\Phi(m_n)\,dx+C\,\sum_i \psi(r_i^{-1})\,r_i^{d-1}\\
&+C\,\sum_i \psi(\|L\|_\infty\,r^{-d}_i)\,\int_{B(x_i,3r_i)}
|db_n|\,dx<\infty, 
\end{split}\]
provided $\phi\leq \tilde \phi$ and recalling that
\[
\sum_i r_i^{d-1}<\infty,\quad \sum_i \int_{B(x_i,3r_i)}
|db_n|\,dx\leq C\int_\Omega |db_n|\,dx.
\]
Denote $\omega_r$ the set of $x$ s.t. $\exists t\in[0,\ T]$ with
$X_n(t,x)\in \Omega_r$. Of course
\[
|\omega_r|\leq \sum_i |\{x,\ \exists t\ X_n(t,x)\in B(x_i,4r_i)\}|
\]
and if $X_n(t,x)\in B(x_i,3r_i)$ then $X_n(s,x)\in B(x_i,5r_i)$ for
$s\in [t-r_i,\ t+r_i]$ and so
\[\begin{split}
|\{x,\ \exists t\ X_n(t,x)\in B(x_i,4r_i)\}|&\leq \frac{1}{r_i}
\int_\Omega \int_0^T \ind_{X_n(t,x)\in B(x_i,5r_i)}\,dt\,dx\\
&\leq \frac{C}{r_i} \int_0^T \int_\Omega \ind_{x\in
  B(x_i,5r_i)}\,dx\,dt\leq C\,r_i^{d-1}.
\end{split}\]
 Hence one has
\begin{equation}
|\omega_r|\leq \eta/8.\label{omegar}
\end{equation}
Because of \eqref{sigma} and Lemma \ref{lemma2} as long as $x$
and $y$ are on the same side of $H$, not in $\omega_r$
 and both at distance
larger than $\delta$ then
\begin{equation}\begin{split}
|b_n(X_n(t,x))-b_n(X_n(t,y))|&=|\tilde b_n(X_n(t,x))-\tilde
b_n(X_n(t,y))|\\
&\leq C \int_{B_K(x,y)}
\tilde m_n\,\left(\frac{1}{|x-z|^{d-1}}+ 
\frac{1}{|y-z|^{d-1}}\right)\,dz,\end{split}\label{estimate}
\end{equation}
with 
\[
\int_\Omega \Phi(\tilde m_n)\,dx<\infty.
\]
\bigskip

\noindent {\em Step 2~: Decomposition of the trajectories.}\\
Fix $w\in S^{d-1}$.
For any $x$ we decompose the time interval $[0,\ T]$ into segments
$]s_i,\ t_i[$ such that on such an interval $X_n(t,x)$ and
$X_{n,\delta}=X_n(t,x+\delta\,w)$ are both on the same side of $H$
and both at a distance larger than $\delta$ of $H$.

On the contrary on each interval $[t_i,\ s_{i+1}]$,  either $X_n(t,x)$
and $X_{n,\delta}=X_n(t,x+\delta\,w)$ 
are on different sides of $H$ or one of them is at
a distance less than $\delta$ of $H$.

Notice that of course $t_i$ and $s_i$ depend on $x$, $\delta$, $b$ and
$\eta$ (as $H$ depends on $\eta$) and 
\[
[0,\ T]=\bigcup_{i\leq n}\ ]s_i,\ t_i[\cup [t_i,\ s_{i+1}].
\]
The transition between the two intervals (at $t_i$ as well as at
$s_i$) is when one of $X_n(t,x)$ or
$X_n(t,x+\delta w)$ is exactly at distance $\delta$ of $H$, because for
the two of them to be either on both side or on the same side of $H$,
one of then has to cross $H$ and thus to come first at a distance $\delta$.  

The idea of the rest of the proof is to bound $|X-X_\delta|$ on
$]s_i,\ t_i[$ with arguments similar to the $W^{1,1}$ case. On $[t_i,\
s_{i+1}]$, the bound will be obtained simply by controlling
$s_{i+1}-t_i$. This will be enough provided that the number of such
intervals $n$ (depending on $x$ again) is not too large.

\bigskip

\noindent {\em Step 3~: Bound on the number of intervals.}\\
At $s_i$ and $t_i$ either $X_n(t,x)$ or $X_{n,\delta}=X(t,x+\delta\,w)$ is at
distance $\delta$ of $H$. Therefore as $\|b\|_\infty\leq 2$, on the
intervals $[s_i-\delta,\ s_i+\delta]$ and $[t_i-\delta,\ t_i+\delta]$,
again either $X_n$ or $X_{n,\delta}$ is at distance less than $3\delta$ of
$H$.

So if there are $n$ such intervals for $x$ then either $X_n$ stays at
distance less than $3\delta$ of $H$ for a total time of at least
$n\delta$ or $X_{n,\delta}$ does.

Now denote by $\omega^1$ the set of all $x$ such that $X_n(t,x)$ stays
distance less than $3\delta$ of $H$ for a total time of at least 
$n\delta$. Obviously
\[
\int_0^T \int_{\omega^1} \ind_{d(X_n(t,x),H)\leq 3\delta}\,dx\,dt\geq
n\delta\,|\omega^1|.
\]
But on the other hand, changing variable
\[\begin{split}
\int_0^T \int_{\omega^1} \ind_{d(X_n(t,x),H)\leq 3\delta}&\,dx\,dt\leq
\int_0^T \int_{\Omega} \ind_{d(X_n(t,x),H)\leq 3\delta}\,dx\,dt\\
\leq &\int_0^T \int \ind_{d(x,H)\leq 3\delta}\,dx\,dt\leq
K\,T\,\delta\,{\cal H}^{d-1}(H)<K\,T\,\delta,
\end{split}\]
as the mesure of $H$ is finite but depends on $\eta$. Thus one finds
\[
|\omega^1|\leq K/n.
\] 
Choose $n=1/(8K\eta)$ so that except a set of measure less than
$\eta/4$ (and equal to $\omega^1\cup(\omega^1-\delta w)$), no trajectory
is decomposed on more than $n$ intervals.

\bigskip

\noindent{\em Step 4~: Control on $[t_i,\ s_{i+1}]$, definitions.}\\
First of all, take any $z\in H$. Denote by $b_-(z)$ the trace of
$b$ on one side of $H$ and $b_+(z)$ on the other side (again
orientation is chosen locally and does not have to hold globally for
$H$). For simplicity, we assume that $b_{-}(z)\cdot\nu(z)\geq 0$ (with
$\nu$ the normal to $H$ chosen according to the orientation).

Now fix some $ \alpha$ to be chosen later.

As $b$ is approximately continuous (see \cite{AFP})  at $z$ on either
side of $H$ one has
\[
r^{-d}\,|\{x\in B_{\pm}(z,r),\ |b(x)-b_{\pm}(z)|\geq
\alpha\}|\longrightarrow 0,
\]
where $B_{\pm}(z,r)$ is the subset of $B(z,r)$ which is on the $-$ or
$+$ side of $H$. 

It is even possible to be more precise. Take a function
$\tilde \eta(r)\rightarrow 0$ as $r\rightarrow 0$ with $\tilde
\eta\leq\eta/(6{\cal H}^{d-1}(H))$. Denote for $r\leq 1/2$
\[
\mu(z,r)_{\pm}=\{x\in B_{\pm}(z,r),\ |b(x)-b_{\pm}(z)|\geq
\alpha\}|,
\]
and
\[
\tilde \mu(z,r)_{\pm}=\{x\in \mu(z,r)_{\pm},\ d(x,H)\geq r\,\tilde\eta/K\},
\]
with $K$ large enough with respect to the Lipschitz constant of $H$.
Assume that $|\mu(z,r)_-|\geq r^d\,\tilde\eta$, then necessarily 
$|\tilde\mu(z,r)_-|\geq r^d\,\tilde\eta/2$; For any
$x\in\tilde\mu(z,r)_-$, one has
\[
\alpha\leq|b(x)-b_-(z)|\leq\int_0^1 |db(\theta x+(1-\theta)z|\,d\theta.
\]
Integrate this inequality over $\tilde \mu$ to find
\[
\int_0^1\int_{B(z,r)\ s.t.\ d(x,H)\geq r\,\tilde\eta/K} 
|db(\theta x+(1-\theta)z)|
\,dx\,d\theta\geq C\,\alpha\,r^d\,\tilde\eta/2.
\]
Now denote $H_r$ the set of $z$ such that   $|\mu(z,r)_{\pm}|\geq
r^d\,\tilde\eta$.  
Integrate the last inequality on $H_r$ on the right hand side and on
$H$ on the left hand side to obtain
\[
\int_H\int_0^1\int_{B(z,r)\ s.t.\ d(x,H)\geq r\,\tilde\eta/K} 
|db(\theta x+(1-\theta)z)|
\,dx\,d\theta\,dz\geq C\,\alpha\,r^d\,\tilde\eta {\cal H}^{d-1}(H_r)/2.
\]
Let us bound the integral on the left hand side. As the case
$\theta\geq 1/2$ is easy, assume $\theta<1/2$ and change variables
locally around $H$ such that $H$ has equation $x_1=0$. The integral is
then 
dominated by
\[\begin{split}
K\,\int_{z\in H}\int_0^{1/2}&\int_{B(z,r), x_1<-r\tilde\eta/K} |db(\theta
x_1, \theta x'+(1-\theta)z)|\,dx\,d\theta\,dz\\
&\leq 2^{d-1}\,K\,\int_{z\in \tilde H}\int_0^{1/2}\int_{B(z,r),
  x_1<-r\tilde\eta/K}|db(\theta
x_1, z)|\,dx\,d\theta\,dz ,
\end{split}\]
with $\tilde H=\{\theta x'+(1-\theta)z,\ z\in H,\ x\in B(z,r)\}$ and
denoting $x=(x_1,\;x')$. Then changing variable from $\theta$ to $\theta
x_1$, one may finally bound by
\[\begin{split}
2^{d-1}\,K\,\int_{z\in \tilde H}\int_{-r/2}^0 &|db(\theta,z)|\,
\int_{-r<x_1<-r\tilde \eta} \frac{dx_1}{\tilde \eta
  r}\,d\theta\,dz\\
&\leq \frac{K}{\tilde \eta}\int_{z\in \tilde
  H}\int_{-r/2}^0 |db(\theta,z)|d\theta\,dz.
\end{split}\] 
Note that this last integral converges to $0$ as $r\rightarrow 0$ and
hence conclude that by choosing $\tilde\eta$ large enough in terms of
$r$, one can ensure that
\[
|H_{r'}|\leq \chi(r'),\quad \chi(r')\rightarrow 0\ \mbox{as}\
r'\rightarrow 0,\quad \chi(r')\leq \eta/C\ \forall r'\leq r.
\]

In addition from the strong convergence $b_n\rightarrow b$, take $n_0$ s.t.
$\forall n\geq n_0$, $\forall z\in H\setminus H_r$ 
\begin{equation}\begin{split}
{r'}^{-d}\,|\{x\in B_{\pm}(z,r'),\ |b_n(x)-b_{\pm}(z)|\geq
\alpha\}|\leq 3\tilde \eta/2\leq &\eta/(4{\cal H}^{d-1}(H))\\
&\quad \forall
|\delta|<r'<r.
\end{split}\label{tildeeta}\end{equation}
This is the second point where $n$ has ot be large enough in terms of $\delta$.

Remark eventually that one may always assume
 that $r\leq \alpha$ (unfortunately not
the opposite way) and $(\alpha\,r)^d<\eta/C$.

\bigskip

\noindent{\em Step 5~: Control on $[t_i,\ s_{i+1}]$, non-crossing
  trajectories.}\\ 
We start by taking apart the trajectories which do not cross $H$. So
define
\[\begin{split}
&H_0=\{z_i\not\in H_r,\ b_{-}(z_i)\cdot\nu(z_i)\leq 2\alpha\},\quad \mbox{and}\\
&\omega^2=\left\{x,\ \exists t\in [0,\ T]\ s.t.\ 
X_n(t,x)\in\bigcup_{z\in H_0} B_{-}(z,\alpha\,r)\cup
B_{+}(z,\alpha\,r)
\right\},\\
&\omega^3=\{x,\ \forall t\in[0,\ T]\ s.t.\ d(X_n(t,x),H)<\alpha\,r/2,\\
&\qquad\qquad\qquad\qquad\qquad
\mbox{one has}\ B(X_n,\alpha\,r)\cap (H\setminus H_r)=\emptyset\}.
\end{split}\]
Note that by \eqref{incompn} (see \ref{incompress} in the introduction), if
$b_{-}(z)\cdot \nu(z)\leq 2\alpha$ on a neighborhood of the border then
$b_{+}\cdot \nu(z)\leq K\alpha$ (with $K$ function of the constant in
\eqref{incomp} and the Lipschitz bound on $H$). Therefore in the
previous definition of $H_0$, one only needs to put $b_{-}$ (in any
case the proof could work the same by essentially ignoring what occurs
on the side where $b\cdot \nu\geq 0$, except in a band of width $\delta$).

Bound $\omega^3$ first. Simply notice that if
\[
\Omega^3=\{x\ \mbox{with}\ d(x,H)<\alpha\,r/2,\
s.t.\ B(x,\alpha\,r)\cap (H\setminus H_r)=\emptyset \} ,
\]
then for $x\in \Omega^3$, ${\cal H}^{d-1}(B(x,\alpha r)\cap H_r)\geq
C\,\alpha^{d-1}\,r^{d-1}$ and as ${\cal H}^{d-1} (H_r)\leq \eta/C$ then 
\[
|\Omega^3|\leq \alpha\,r\eta/C,\quad |\Omega^3_r|=
|\{x,\ d(x,\Omega^3)\leq \alpha r\}|\leq \alpha\,r\,\eta/C.
\]
As $\omega^3=\{x,\ \exists t\in [0,\ T]\ X_n(t,x)\in \Omega^3\}$, and
if $X_n(t,x)\in\Omega^3$ then $X_n(s,x)\in \Omega^3_r$ for $s\in [t,\
t+\alpha r]$
\[
\alpha\,r|\omega^3|\leq \int_0^T \int_{\Omega} \ind_{X_n(t,x)\in
  \Omega^3_r}\,dx\leq C\int_0^T\int_{\R^d} \ind_{x\in
  \Omega^3_r}\,dx\leq T\,\alpha\,r\,\eta/C,
\]
by change of variable and the bound on $|\Omega^3_r|$. This implies
\begin{equation}
|\omega^3|\leq T\,\eta/C.\label{omega3}
\end{equation}

Let us now bound $\omega^2$, decompose
\[
\omega^2=\omega^{2,1}\cup \omega^{2,2},
\]
with $\omega^{2,1}$ the subset of all $x\in \omega^2$ s.t. the
trajectory $X_n(t,x)$ stays at least a time interval $r/2$ 
at a distance less than $K\alpha\,r$ of $H$. Notice that
\[
\int_{\omega^{2,1}} \int_0^T \ind_{d(X_n,H)\leq K\alpha\,r}\,dt\,dx\geq
|\omega^{2,1}|\,\frac{r}{2}.
\]
On the other hand by \eqref{incompn}
\[\begin{split}
\int_{\omega^{2,1}} \int_0^T \ind_{d(X_n(t,x),H)\leq K\alpha\,r}\,dt\,dx&\leq
\int_0^T \int_{\R^d} \ind_{d(x,H)\leq K\alpha\,r}
\,dx\,dt\\
&\leq K\,T\,\alpha\,r\,{\cal
  H}^{d-1}(H).
\end{split}\]
Therefore 
\[
|\omega^{2,1}|\leq \eta/8,
\]
provided that
\begin{equation}
C\,T\,K\,{\cal H}^{d-1}(H)\,\alpha\leq \eta/8.
\end{equation}

There remains $\omega^{2,2}$ which is made of those $x\in\omega^2$
staying a time less than $r/2$ at a distance less than $K\alpha r$
of $H$. Denote $t_0$ s.t. $X_n(t_0,x)\in B_{-}(z,\alpha\,r)$ for some
$z$ (and of course the same analysis holds for the $+$ side).

If this is so then the average over the interval $[t_0,\ t_0+r/2]$ 
\[
\frac{2}{r}\int_{t_0}^{t_0+r/2} \nu(X_n)\cdot b_n(X_n(t,x))\,dt\geq
(K-1)\,\alpha. 
\]
Taking $K$ large enough with respect to the Lipschitz bound on $H$,
this implies
\[
\frac{2}{r}\int_{t_0}^{t_0+r/2} \nu(z)\cdot b_n(X_n(t,x))\,dt\geq
K\,\alpha, 
\]
and consequently
\[
\frac{2}{r}\int_{t_0}^{t_0+r/2}  |b_n(X_n(t,x))-b_{-}(z)|\,dt\geq
K\,\alpha.
\]
As $1\leq b_1\leq 2$, notice that on the time interval $[t_0,\ t_0+r/2]$ then  
$X_n(t,x)$ stays inside the ball $B(z,r)$ and moreover the length of
the path is of order $r$. Consequently, denoting $\omega^{2,2}_z$ the
subset of $\omega^{2,2}$ corresponding to the ball $B(z,r)$, one has
by change of variable (and again the use of \eqref{incompn})
\[
\frac{1}{r^d}\{x\in B_{-}(x),\ |b_n(x)-b_{-}(z)|\geq \alpha\}\geq
K\,|\omega^{2,2}_z|,
\]
which finally gives $\omega^{2,2}_z\leq \eta r^d/K$ and summing over
all $z=z_i$ 
\begin{equation}
|\omega^2|\leq |\omega^{2,1}|+|\omega^{2,2}|\leq \eta/8+\eta\,\sum_i
r^d/K\leq \eta/4,\label{omega2} 
\end{equation}
by choosing $K$ large enough.

\bigskip

\noindent{\em Step 6~: Control on $[t_i,\ s_{i+1}]$, the crossing
  trajectories.}\\ 
We now only have to take into account the trajectories passing
through. Consider for example
\[
\begin{split}
&H_1=\{z\not\in H_r,\ b_{-}(z)\cdot\nu(z)\geq 2\alpha\ 
\},\quad \mbox{and}\\
&\omega^0=\left\{x,\ \exists t\in [0,\ T]\ s.t.\ 
X_n(t,x)\in\bigcup_{z\in H_1} B_{-}(z,\alpha\,r)
\right\}.\end{split}
\]
The same holds if one uses
$b_+(z_i)\cdot\nu(z)<-2\alpha$ in the definition 
and note again that by \eqref{incomp}
one has necessarily one or the other as the case
$|b_{\pm}\cdot\nu(z)|\leq2\alpha$ was taken care of in the previous
step. 

As $|H_{r'}|\rightarrow 0$, there exists an extracted sequence 
$r_k\rightarrow 0$ (with
$r_0=r$) s.t. for any $0<\gamma<1$
\begin{equation}\begin{split}
&\sum_k r_k<\infty,\quad
\sum_k |H_{r_k}|<\infty,\quad\mbox{and so}\ \sum_{k=0}^\infty
|H_{r_k}|\leq \eta/K,\\
& \sum_k \tilde \eta(r_k)<\infty,\qquad \sum_{k=0}^{+\infty} \tilde
\eta(r_k)\leq \eta/K.
\end{split}\label{rk1}\end{equation}
For any $\delta$, denote $k_\delta$ s.t.
\begin{equation}
\sum_{k\leq k_\delta} \frac{\delta}{r_k}\leq \eta/K. \label{rk2}
\end{equation}
Of course $r_{k_\delta}>\delta$ but note that
$k_\delta\rightarrow+\infty$ as $\delta\rightarrow 0$.

\medskip 

Denote $\Omega_+$ the set of
$y\in \bigcup_{z\in H}B_+(z,r)$ s.t. $d(y,H)>\delta$ and
$\Omega_-=\bigcup_{z\in H}B_+(z,r)\setminus \Omega_+$.

For any scale $r_k>C\delta$ with $k\geq 1$,
define 
$\omega_{k}$ as the set of $x$ s.t. \\
$(i)$ $\exists t_0\in[0,\ T]$ with $d(X_n(t_0,x),H)\leq r_k\alpha/K$,
$K\geq 6$,\\
$(ii)$ one never has $X_n(t,x)\in \Omega_+$ for any $t\in [t_0,\
t_0+r_k/3]$.

By $(i)$ and $\|b_n\|_\infty\leq 2$, one has that $X_n(t,x)$ remains in
the same ball  $B(X_n(t_0,x),r_k/2)$
for all $t\in[t_0,\ t_0+r_k/K]$ (for $K\geq 4$). The intersection of this ball 
with $H$ has diameter larger than $r_k/3$ as $d(X_n(t_0,x),H)\leq
r_k\alpha/6$.

Accordingly decompose $\omega_{k}=\omega_k^1\cup\omega_k^2$ with
$\omega_k^1$ the set of $x\in\omega_k$ s.t. 
\[
B(X_n(t_0,x),r_k/2)\cap (H\setminus H_{r_k})=\emptyset,
\]
and $\omega_k^2=\omega_k\setminus \omega_k^1$.

Start by bounding $\omega_k^1$. Define
\[
\Omega_k^1=\{x,\ B(x,r_k/2)\cap (H\setminus H_{r_k})=\emptyset\}.
\]
Just as in the previous step
\[
|\Omega_k^1|\leq C\,{\cal H}^{d-1}(H_{r_k})\,r_k,\quad |\tilde \Omega_k^1|=\{x,\
d(x,\Omega_k^1)\leq r_k\}|\leq C\,{\cal H}^{d-1}(H_{r_k})\,r_k.
\]
On the other hand $X_n(t,x)$ stays inside $\tilde\Omega_k^1$ for all
the interval $[t_0,\ t_0+r_k/K]$ so again
\[
|\omega_k^1|\,r_k/K\leq \int_0^T \int_\Omega \ind_{X_n(t,x)\in \tilde
  \Omega_k^1} \,dx\,dt\leq \int_0^T\int_{\R^d} \ind_{x\in\tilde
  \Omega_k^1}\,dx\,dt\leq C\,{\cal H}^{d-1}(H_{r_k})\,r_k,
\]
which gives the desired bound
\begin{equation}
|\omega_k^1|\leq K\,{\cal H}^{d-1}(H_{r_k}).\label{omegak1}
\end{equation}

\medskip

Now for $\omega_k^2$, notice that for any $x\in \omega_k^2$ there
exists $z\in B(X_n(t_0,x),r_k/2)\cap (H\setminus H_{r_k})$. Hence we
may take $N$ points $z_i\in H\setminus H_{r_k}$ with
\[
N\leq K\,{\cal H}^{d-1}(H)/r_k^{d-1},
\] 
and such that for any $x\in \omega_k^2$, there is $z_i\in
B(X_n(t_0,x),r_k/K)$. 

On the interval $[t_0,\ t_0+r_k/K]$, $X_n$ is in $\Omega_-$. Denote
$H^\delta=\{x,\ d(x,H)\leq \delta\}$, decompose again $\omega_k^2$
into $\omega_k^{2,\delta}\cup\tilde \omega_k^2$ with
$\omega_k^{2,\delta}$ the set of $x$ such that $X_n$ stays more than
a total time $\Delta=\alpha r_k/K$ in $H^\delta$. Bound
\[
|\omega_k^{2,\delta}|\,\alpha r_k/K\leq \int_0^T\int_\Omega
\ind_{d(X_n(t,x),H)<\delta}\,dx\,dt\leq K\,\delta\,{\cal H}^{d-1}(H).
\]
Thus
\begin{equation}
|\omega_k^{2,\delta}|\leq K\,\delta\,{\cal H}^{d-1}(H)/r_k.\label{omegak21}
\end{equation}
For $x\in\tilde\omega_k^2$, denote by $I$ the subset of $[t_0,\ r_k/K]$
s.t. $X_n(t,x)\not\in H^\delta$. As $X_n$ never reaches $\Omega_+$ and
stays in $H^\delta$ at most $\alpha r_k/K$ one has
\[
\frac{1}{|I|}\int_I b_n
(X_n(t,x))\cdot \nu(X_n)\,dt< \frac{K}{r_k}\;(\delta+\alpha
r_k/K+2\,\alpha r_k/K)< \alpha/3.
\]
Moreover for any $t\in[t_0,\ t_0+r_k/K]$, one has $|X_n(t,x)-z_i|\leq
r_k/K$ and as $\nu$ is lipschitz (at least around $H$) then
\[
\frac{1}{|I|}\int_I b_n(X_n(t,x))\cdot \nu(z_i)\,dt<\alpha/2,
\]
and
\[
\frac{1}{|I|}\int_I \ind_{|b_n(X_n(t,x))-b_\pm(z_i)|\geq \alpha}\,dt>1.
\]
Put $\tilde \omega_k^i=\{x\in \tilde\omega_k^2,\ X(t_0,x)\in B(z_i,r_k/K)$.
By \eqref{tildeeta} and integrating along the trajectories (as we did many
times before), we deduce that
\[
| \tilde\omega_k^i|\leq K\tilde \eta(r_k)\,r_k^{d-1}.
\]
Summing over $i$,
\begin{equation}
|\tilde\omega_k^2|\leq K\tilde \eta(r_k).\label{omegak2}
\end{equation}

\bigskip

\noindent{\em Summary and conclusion of the estimate.}\\
Define
\[
\tilde \omega=\omega_r\cup\omega^1\cup \omega^2\cup\omega^3\cup_k
\omega_k,\quad \bar\omega=\tilde\omega\cup (\tilde\omega-\delta w).
\]
By the previous computations (end of steps 1, 3, step 5 \eqref{omega3}
and \eqref{omega2}
and step 6 \eqref{omegak1}, \eqref{omegak21}, \eqref{omegak2})
\begin{equation}
|\bar\omega|\leq \eta/4+\eta/4+T\eta/C+\eta/4+K\,\sum_k (|H_{r_k}|+\delta/r_k+
\tilde \eta(r_k))\leq 7\eta/8,\label{omega}
\end{equation}
by \eqref{rk1} and \eqref{rk2}.

Now for $x\not\in \bar\omega$, we have for any $t\in [s_i,\ t_i]$
\[
|b(X_n)-b(X_{n,\delta})|\leq C\int_{B_K(X_n,\ X_{n,\delta})} \tilde m_n\left(
\frac{1}{|X_n-z|^{d-1}}+\frac{1}{|X_{n,\delta}-z|^{d-1}}\right)\,dz.
\]

As in the proof of \ref{W11}, we define
\[
\psi(|\delta|)=\inf_M M-2K\frac{M}{\phi(M)}\,\log(|\delta|).
\]
And we obtain
\[\begin{split}
\int_{\bar\omega^c}&\sum_i
\log\left(\frac{\delta+|X_n(t_i,x)-X_{n,\delta}(t_i,x)|}
{\delta+|X_n(s_i,x)-X_{n,\delta}(s_i,x)|}\right)\\
&\quad \leq \int_{\omega^c}\sum_i \int_{s_i}^{t_i}
\frac{|b(X_n)-b(X_{n,\delta})|}{|\delta|+|X_n-X_{n,\delta}|}\,dt\,dx\\
&\leq
C\int_\Omega\int_0^T \int_{B_K(X_n,\ X_{n,\delta})} \tilde m_n\left(
\frac{1}{|X_n-z|^{d-1}}+\frac{1}{|X_{n,\delta}-z|^{d-1}}\right)\,dz\leq
\psi(|\delta|). 
\end{split}\]
So up to the first time $t_1$, one has that
\[
|X_n-X_{n,\delta}|\leq \delta\,\exp(8\psi(|\delta|)\,n/\eta),
\]
except for $x\in \omega^4_1\cup\bar\omega$ with
\[
|\omega^4_1|\leq \eta/8n.
\]
Now by induction let us prove that up to time $t_i$
\[
|X_n-X_{n,\delta}|\leq \delta_i,
\]
except for $x\in \bar\omega\cup\omega^4_i$ with
\[
|\omega_i^4|\leq i\,\eta/8n,
\]
and $\delta_i$ defined through
\begin{equation}
\delta_{i+1}=(\delta_i+r_{k_i}/K)\,\exp(8n\,
\psi(\delta_i+r_{k_i}/K)/\eta),\quad r_{k_i}=\inf \{r_k,\
\alpha\,r_k/K>\delta+\delta_i\}.
\end{equation}
This is true for $i=1$. So assume it is still true up to $t_i$ and
study what happens on $[t_i,\ t_{i+1}]$. First of all on $[t_i,
s_{i+1}]$~: at $t_i$ assume for instance that
$d(X_{n,\delta},H)=\delta$ and take $k_i$ s.t. $r_{k_i}=\inf\{r_k,
\alpha r_k/K>\delta+\delta_i\}$. Then necessarily $d(X_n(t_i),H)\leq
\alpha r_{k_i}/K$. As $x\not\in \bar\omega$, by step 6, one has that
$s_{i+1}\leq t_{i}+r_{k_i}/K$. 
Therefore 
\[
|X_n(s_i)-X_{n,\delta}(s_i)|\leq \delta_i+r_{k_i}/K.
\]
Now on $[s_{i+1},\ t_{i+1}]$, simply write
\[\begin{split}
&\int_{\bar\omega^c\cap \omega_i^4}
\log\left(\frac{\delta_i+r_{k_i}/K+|X_n(t_{i+1},x)-X_{n,\delta}(t_{i+1},x)|}
{\delta_i+r_{k_i}/K+|X_n(s_{i+1},x)-X_{n,\delta}(s_{i+1},x)|}\right)\\
&\quad \leq \int_{\omega^c\cap\omega_i^4} \int_{s_{i+1}}^{t_{i+1}}
\frac{|b(X_n)-b(X_{n,\delta})|}{\delta_i+r_{k_i}/K+|X_n-X_{n,\delta}|}\,dt\,dx\\
&\leq
C\int_\Omega\int_0^T \int_{B_K(X_n,\ X_{n,\delta})} \tilde m_n\left(
\frac{1}{|X_n-z|^{d-1}}+\frac{1}{|X_{n,\delta}-z|^{d-1}}\right)\,dz\\
&\qquad\leq
\psi(\delta_i+r_{k_i}/K). 
\end{split}\]
Therefore for $t\leq t_{i+1}$, one has that
\[
|X_n-X_{n,\delta}|\leq \delta_{i+1}= (\delta_i+r_{k_i}/K)\,\exp(8n\,
\psi(\delta_i+r_{k_i}/K)/\eta),
\]
except for $x\in\bar\omega\cup\omega_{i+1}^4$, with $\omega_i^4\subset
\omega_{i+1}^4$ and
\[
|\omega_{i+1}^4\setminus\omega_i^4|\leq \eta/8n,\ s.t.\
|\omega_{i+1}^4|\leq (i+1)\,\eta/8n.
\]
Finally for any $t\in [0, T]$, one has that
\[
|X_n-X_{n,\delta}|\leq \delta_n,
\]
for any $x\not\in \omega$ with $\omega=\bar\omega\cup\omega_n$ and
thus
\[
|\omega|\leq \eta.
\]
This concludes the proof once one notices that $\delta_n\rightarrow 0$
as $\delta\rightarrow 0$ since $\psi(\delta)/|\log \delta|\rightarrow
0$.

\bigskip

{\em Compactness.}\\
 Let us just briefly indicate how to obtain the
compactness from the estimate on $X_n-X_{n,\delta}$.

Notice first that in the previous proof, $n$ had to be taken large
enough in terms of $\delta$. So if for a given $n$, one considers
$\delta'<<\delta$ 
instead of $\delta$, it is not true that
$|X_n-X_{n,\delta'}|$ may be bounded by $\ep(\delta')$. Instead the best
that can be done is nothing as long as $|X_n-X_{n,\delta'}|\leq
\delta$ and bound as before once $|X_n-X_{n,\delta'}|\geq \delta$. So
in the end $|X_n-X_{n,\delta'}|$ is only bounded by $\ep(\delta)$.

For the compactness of $X_n$, recalling that $X_n$ is uniformly
lipschitz in time, we apply the usual criterion saying that
$X_n$ is compact in $L^1([0,\ T]\times\Omega)$ iff
\[
\forall \gamma,\, \forall w,\ \exists \delta\ s.t.\ \forall
\delta'<\delta,\quad
\sup_n \int_0^T\int_\Omega |X_n(t,x)-X_n(t,x+\delta'w)|\,dx\,dt<\gamma. 
\]
So fix $\gamma$ and $w$. First choose $\eta$ s.t.
\[
\eta\,T\,\max_\Omega |x|<\gamma/4.
\]
Now apply the previous quantitative estimate to obtain a function
$\ep(\delta)$ (recall that this depends on $\eta$). Choose $\delta_1$
with $\ep(\delta_1)<\gamma/(2|\Omega|T)$. This gives $n_1$ s.t. for
any $n\geq n_1$ and any $\delta<\delta_1$
\[
|X_n-X_{n,\delta}|\leq \ep(\delta_1)<\gamma/(2|\Omega|T),\ \forall
x\in \Omega\setminus\omega_\delta\ \mbox{with}\ |\omega_\delta|\leq \eta.
\]
Now for $n<n_1$ as $X_n$ is regular (not uniformly in $n$ but there
are only $n_1$ indices $n$ now), choose $\delta_2$ s.t. for any
$n<n_1$ and any $\delta<\delta_2$
\[
|X_n-X_{n,\delta}|<\gamma/(|\Omega|T).
\]
Consequently take $\delta=\min(\delta_1,\delta_2)$. For any $n$ and
any $\delta'<\delta$, if $n<n_1$ then simply
\[
\int_0^T \int_\Omega |X_n(t,x)-X_{n}(t,x+\delta'w)|\,dx\,dt\leq
T\,|\Omega|\,\sup_x |X_n-X_{n,\delta}|<\gamma.
\]
And if $n\geq n_1$, then decompose
\[\begin{split}
\int_0^T \int_\Omega &
|X_n(t,x)-X_{n}(t,x+\delta'w)|\,dx\,dt=\int_0^T\int_{\omega_{\delta'}}\ldots 
+\int_0^T\int_{\Omega\setminus\omega_{\delta'}}\ldots\\
&\leq 2\,T\,|\omega_{\delta'}|\,\max_\Omega |x|
+T\,|\Omega|\,\sup_{\Omega\setminus\omega_{\delta'}} |X_n-X_{n,\delta'}|<\gamma. 
\end{split}\]
Hence the compactness criterion is indeed satisfied.
\section{The $BV$ case}
I do not know how to perform a rigourous proof in the full $BV$ case
so the purpose of this section is only to try to explain what can be
done and where are the problems. Accordingly most technical details
are omitted. 

The aim here is the uniqueness of the flow, which turns out
to be much simpler than the compactness (as in the $SBV$ case). So
consider two solutions $X(t,x)$ and $Y(t,x)$ to \eqref{ode}, both of
them satisfying \eqref{incomp} and \eqref{assumpb}. We can define the
set of points where $X$ and $Y$ start being different and as both are
flows (semi-groups) then it comes to
\[
F=\{x\;|\ \exists t_n\rightarrow 0\ s.t.\ X(t_n,x)\neq Y(t_n,x)\}.
\]
If ${\cal H}^{d-1} (F)=0$ then everything is fine as almost no
trajectory passes through $F$ and uniqueness holds for $a.e.$ initial
data $x$. So we may assume that  ${\cal H}^{d-1} (F)>0$.

Note first that the Lebesgue measure of $F$ is necessarily $0$ because
non uniqueness may occur only where $db\not\in L^1$ and this set has
vanishing Lebesgue measure. 

On an interval $[t_0,\ s_0]$, $|X(t,x)-Y(t,x)|$ passes from
$\delta$ to $2\delta$ (the infimum is $\delta$ and the maximum larger than
$2\delta$) then by Lemma \ref{lemma} 
\[
\int_{t_0}^{s_0} \int_{B(X,2\delta)} |db(z)|\geq 1/C.
\]
Therefore defining 
\[
F_\ep=\{x,\ d(x,F)\leq \ep\},
\]
one has for any $\ep>0$ that
\[
\int_{F_\ep} |db(z)|\geq 1/C.
\]
It implies that some mass of $db$ is concentrated on $\bar F$ (the
closure of $F$) or $\int_{\bar F} |db(z)|>0$. As $|db|$ is the
distributional derivative of a $BV$ function, it cannot concentrate
mass on a purely unrectifiable set. Therefore if $\bar F$ has a purely
unrectifiable set $F_0$ then almost no trajectory $X$ or $Y$ crosses
$F_0$. So we may reduce ourselves to the case where $\bar F$ does not
contain any such set (by considering $\bar F\setminus F_0$).

Take $H\subset \bar F$, rectifiable and with $0<{\cal
  H}^{d-1}(H)<\infty$. Denote $\nu(x)$ the normal on $H$ and look at
the part of $H$ where the trajectories $X$ and $Y$ are not tangent
\[
H_0=\{x\in H,\ b(x)\cdot \nu(x)\neq 0\}.
\]
If ${\cal H}^{d-1}(H_0)> 0$ then necessarily the set 
\[\Omega_H=\{x,\
\exists t\ X(t,x)\in H\ or\ Y(t,x)\in H\}\]
has non zero Lebesgue
measure. Consequently, by the same arguments as before
\[
\int_H |db(z)|>0,
\]
and this implies that ${\cal H}^{d-1}(J\cap H)>0$ where $J$ is the
jump set of $b$. This case was dealt with before as it is exactly the
$SBV$ situation.

So considering $\tilde F=\bar F\setminus J$ instead of $\bar F$, the
only remaining situation is where for any rectifiable $H\subset \tilde
F$ with ${\cal H}^{d-1}(H)<\infty$ then 
\[
{\cal H}^{d-1}\left(\{x\in H,\ b(x)\cdot \nu(x)\neq 0\}\right)=0,
\]
or for ${\cal H}^{d-1}$ $x$ in $\tilde F$, one has $\nu(x)\cdot b(x)=0$ but
still $\Omega_{\tilde F}=\{x,\ \exists t\ X(t,x)\in \tilde F\ or\
Y(t,x)\in \tilde F\}$ has non zero Lebesgue measure. This is the case
which cannot be handled. Note that, rather unsurprisingly, 
the structure of the problem here
is very similar to the one faced in \cite{ADM}.
 
\bigskip

{\bf Acknowledgements.}\\
I am much indebted to C. DeLellis for having introduced the problem to
me. Quite a few ideas in section \ref{sec2d} were inspired by
an attempt to solve another Bressan's conjecture
 with C. DeLellis  and U. Stefanelli. Finally I wish to thank
 S. Bianchini, F. Bouchut, Y. Brenier and particularly L. Ambrosio
for fruitful comments.



\begin{thebibliography}{XX}
\bibitem{Ai} M. Aizenman, On vector fields as generators of 
flows: A counterexample to Nelson's conjecture. {\em Ann. Math.} (2)
{\bf 107}
(1978), pp. 287--296.  

\bibitem{Am} L. Ambrosio, Transport equation and Cauchy problem for
  $BV$ vector fields. {\em Invent. Math.} {\bf 158}, 227--260 (2004).

\bibitem {AC} L. Ambrosio, G. Crippa, Existence, uniqueness, stability and
 differentiability properties of the flow associated to weakly
 differentiable vector fields. Lecture notes of the Unione Matematica
 Italiana, Springer Verlag, to appear.

\bibitem{ADM} L. Ambrosio, C. De Lellis, J. Mal\'y, On the chain rule for the
divergence of vector fields: applications, partial results, open
problems,  Perspectives in nonlinear partial differential equations,
31--67, {\em Contemp. Math.}, {\bf 446}, 
Amer. Math. Soc., Providence, RI, 2007.

\bibitem{AFP} L. Ambrosio, N. Fusco, D. Pallara. Functions of Bounded
  Variation and Free Discontinuity Problems.Oxford 2000.
 
\bibitem{ALM} L. Ambrosio, M. Lecumberry, S. Maniglia,
 Lipschitz regularity and approximate differentiability of the
 DiPerna-Lions flow. {\em Rend. Sem. Mat. Univ. Padova} {\bf 114} (2005),
 29--50.

\bibitem{Bo} F. Bouchut, Renormalized solutions to the Vlasov 
equation with coefficients of bounded variation. {\em
 Arch. Ration. Mech. Anal.} {\bf 157} (2001), pp. 75--90.

\bibitem{BC} F. Bouchut, G. Crippa, Uniqueness, renormalization, 
and smooth approximations for linear transport equations.  {\em SIAM
J. Math. Anal.}  {\bf 38}  (2006),  no. 4, 1316--1328. 

\bibitem{BD} F. Bouchut, L. Desvillettes, 
On two-dimensional Hamiltonian transport equations with continuous
coefficients. {\em Diff. Int. Eq.} (8) {\bf 14} (2001), 1015--1024. 

\bibitem{BJ} F. Bouchut, F. James, 
One dimensional transport equation with discontinuous
coefficients. {\em Nonlinear Anal.} {\bf 32} (1998), 891--933. 

\bibitem{Br1} A. Bressan, An ill posed Cauchy problem for 
a hyperbolic system in two space dimensions.
 {\em Rend. Sem. Mat. Univ. Padova} {\bf 110} (2003), 103--117. 

\bibitem{Br2} A. Bressan, A lemma and a conjecture on the cost of 
rearrangements, {\em Rend. Sem. Mat. Univ. Padova} {\bf 110} (2003),
97--102. 

\bibitem{CP} I. Capuzzo Dolcetta, B. Perthame, On some analogy between different
approaches to first order PDE's with nonsmooth coefficients.
{\em Adv. Math. Sci. Appl.} {\bf 6} (1996), 689-- 703. 

\bibitem{CCR} F. Colombini, G. Crippa,  J. Rauch, A note on two-dimensional
transport with bounded divergence. {\em Comm. Partial Differential
  Equations} {\bf 31} (2006), 1109--1115.

\bibitem{CL} F. Colombini, N. Lerner, 
Uniqueness of continuous solutions for BV vector fields. {\em Duke
Math. J.} {\bf 111} (2002), 357--384. 

\bibitem{CL2} F. Colombini, N. Lerner,
Uniqueness of $L\sp \infty$ solutions for a class of conormal $BV$
vector fields. Geometric analysis of PDE and several
complex variables, 133--156, 
{\em Contemp. Math.} {\bf 368}, Amer. Math. Soc., Providence, RI, 2005.

\bibitem{Cr} G. Crippa, The ordinary differential equation with
  non-Lipschitz ve ctor fields.  {\em Boll. Unione Mat. Ital.} (9) {\bf 1}
  (2008),  no. 2, 333--348.  

\bibitem{CD} G. Crippa, C. DeLellis, Estimates and regularity results
  for the DiPerna-Lions flow. 
{\em J. Reine Angew. Math.} {\bf 616} (2008), 15--46.

\bibitem{Da} C. M. Dafermos. Generalized characteristics in hyperbolic
  systems  of conservation laws. {\em Arch. Rat. Mech. Anal.}  {\bf
    107} (1989), 127--155.

\bibitem{DeL} C. De Lellis, Notes on hyperbolic systems of 
conservation laws and transport equations. Handbook of differential
equations, Evolutionary equations, Vol. 3 (2007).
 

\bibitem{DeP} N. De Pauw, Non unicit\'e des solutions born\'ees pour 
un champ de vecteurs $BV$ en dehors d'un
hyperplan. {\em C.R. Math. Sci. Acad. Paris} {\bf 337} (2003), 249--252. 

\bibitem{DL} R.J. DiPerna, P.L. Lions, Ordinary differential
  equations, 
transport theory and Sobolev spaces. {\em Invent. Math.} {\bf 98}
(1989), 511--547. 

\bibitem{Fi} A. F. Filippov. Differential equation with 
discontinous right-hand side. {\em Amer. Math. Soc. Transl.} {\bf 42}
(1962), 199--231.

\bibitem{Ha1} M. Hauray, On two-dimensional Hamiltonian transport 
equations with $\mathbb{L}_{loc}^p$ 
coefficients. {\em Ann. IHP. Anal. Non Lin.} (4) {\bf 20} (2003), 625--644.

\bibitem{Ha2} M. Hauray,
On Liouville transport equation with force field in $BV\sb {\rm loc}$.
{\em Comm. Partial Differential Equations} {\bf  29} (2004), no. 1-2, 207--217. 

\bibitem{HLL} M. Hauray, C. Le Bris, P.L. Lions, Deux remarques 
sur les flots g\'en\'eralis\'es d'\'equations diff\'erentielles
ordinaires. {\em C. R. Math. Acad. Sci. Paris}  {\bf 344}  (2007),
no. 12,
 759--764.

\bibitem{LL} C. Le Bris, P.L. Lions, Renormalized solutions of some
  transport equations with partially $W^{1,1}$ velocities and
  applications. {\em Ann. Mat. Pura Appl.} {\bf 183} (2004), 97--130.

\bibitem{Li} P.L. Lions, Mathematical topics in fluid mechanics, 
Vol. I: incompressible models. Oxford Lect. Ser. Math. Appl. 3
(1996). 

\bibitem{Li2} P.L. Lions, Sur les \'equations diff\'erentielles
  ordinaires 
et les \'equations de transport. {\em C. R. Acad. Sci., Paris, Sér. I,
  Math.}
 {\bf 326} (1998), 833--838.

\bibitem{PP} G. Petrova, B. Popov, Linear transport equation with 
discontinuous coefficients. {\em Comm. Partial Differential Equations}
{\bf 24} (1999),
1849--1873.

\bibitem{PR} F. Poupaud, M. Rascle, Measure solutions to the liner 
multidimensional transport equation with non-smooth
coefficients. {\em Comm. Partial Differential Equations} {\bf 22} (1997), 
337--358.

\end{thebibliography}
\end{document}